\documentclass[3p,12pt]{elsarticle}
\usepackage{amsmath}
\usepackage{cases}
\usepackage{booktabs}
\usepackage{multirow}
\usepackage{xcolor}
\usepackage[normalem]{ulem}
\usepackage{tabularx}
\usepackage{siunitx}
\usepackage[caption=false]{subfig}
\usepackage{comment}
\usepackage[ruled, linesnumbered]{algorithm2e}
\usepackage{algpseudocode}
\usepackage{grffile}
\usepackage{rotating}
\usepackage{lineno}
\usepackage{hyperref}
\usepackage{graphicx,enumerate}
\usepackage{amssymb}
\allowdisplaybreaks[4]
\usepackage{bm}
\usepackage[caption=false]{subfig}
\usepackage{caption,color}
\usepackage{subeqnarray}
\usepackage{multirow}
\usepackage{lscape}
\usepackage{rotating}
\usepackage{graphics}
\usepackage{epstopdf}
\usepackage{threeparttable,booktabs,tabularx}
\journal{Journal of Computational Physics}
\usepackage{floatrow}
\usepackage{floatrow}\newfloatcommand{capbtabbox}{table}[][\FBwidth]
\floatname{algorithm}{Algorithm}
\biboptions{numbers,sort&compress} 

\newcommand{\bq}{\bm q}

\newcommand{\bu}{\bm u}

\newcommand{\bxi}{\bm \xi}

\newcommand{\Kn}{\text{Kn}}

\begin{document}
\begin{frontmatter}
\title{A general synthetic iterative scheme with high-order spatial discretization for rarefied gas mixture flows}

\author[hkust1]{Jianan Zeng}
\author[hkust1,hkust2]{Wei Su\corref{mycorrespondingauthor}}
\cortext[mycorrespondingauthor]{Corresponding author}
\ead{weisu@ust.hk}

\address[hkust1]{Division of Emerging Interdisciplinary Areas, Academy of Interdisciplinary Studies, The Hong Kong University of Science and Technology, Hong Kong SAR, China}
\address[hkust2]{Department of Mathematics, The Hong Kong University of Science and Technology, Hong Kong SAR, China}

\begin{abstract}

A high-order discontinuous Galerkin formulation of the general synthetic iterative scheme (GSIS) is developed for steady rarefied gas mixture flows.
The resulting GSIS-DG discretizes the mesoscopic kinetic equations for monatomic gas mixtures and the macroscopic synthetic equations with the same polynomial degree and solves them alternately.
The kinetic solution supplies nonequilibrium corrections to the stress and heat flux, while the synthetic equations accelerate the convergence toward the steady state. 
Asymptotic analysis shows that GSIS-DG recovers the continuum limit without requiring the spatial mesh to resolve the molecular mean free path.
Numerical results over a wide range of Knudsen numbers demonstrate high-order spatial accuracy and rapid steady-state convergence, with more than an order-of-magnitude reduction in computational cost relative to conventional kinetic iteration in the near-continuum regime.
Comparisons with second-order finite-volume GSIS further demonstrate that GSIS-DG achieves comparable or better accuracy on coarser physical meshes with fewer total spatial degrees of freedom.

\end{abstract}


\end{frontmatter}
\section{Introduction}\label{sec:1}

Gas mixtures are widely encountered in modern engineering problems, such as astronautics~\cite{votta2013hypersonic}, vacuum technologies~\cite{benschop2008extreme,fu2019euv}, and micro-electromechanical systems~\cite{Beskok_book}. These flows often exhibit rarefaction effects, where the Knudsen number, defined as the ratio of the molecular mean free path to the characteristic length, is not negligible. For example, in the dynamic gas lock used in extreme ultraviolet (EUV) lithography, clean gas is continuously injected to reduce the partial pressure of contaminants inside the vacuum optical cavity~\cite{teng2023pollutant}, where the molecular masses and species concentrations of the gas mixture may differ by orders of magnitude. Due to the difficulties of experimental measurements in low-pressure environments, numerical simulation provides an important tool for understanding and optimizing rarefied gas mixture flows. In the presence of significant rarefaction effects, the Navier--Stokes (NS) equations become inadequate. Instead, the Boltzmann equation or simplified kinetic model equations for gas mixtures derived from gas kinetic theory~\cite{McCormack1973,Andries2002JSP,Kosuge2009EJMB,Brull2015CMS,Groppi2011EPL,Todorova2019EJMB,li2024kinetic} should be employed to describe flows spanning the continuum, slip, transition, and free-molecular regimes.

The Boltzmann equation and kinetic model equations can be solved by stochastic or deterministic approaches, among which the direct simulation Monte Carlo (DSMC) method~\cite{bird1994molecular} and the discrete velocity method (DVM)~\cite{broadwell1964discrete} are widely used for rarefied gas flows. However, both approaches become computationally expensive in the near-continuum regime. In DSMC, the spatial cell size and time step are typically constrained by the molecular mean free path and mean collision time, respectively~\cite{bird1994molecular,fei2022usppoly}, whereas deterministic kinetic solvers based on conventional iterations suffer from slow convergence toward the steady state as the Knudsen number decreases. For mixtures with disparate molecular masses, the fast motion and collisions of the light species require small time steps, whereas equilibration between species can take much longer~\cite{li2024kinetic}. These challenges have motivated the development of multiscale methods, including the (discrete) unified gas kinetic scheme (UGKS)~\cite{xu2010unified,guo2013discrete}, the Fokker--Planck DSMC method~\cite{Gorji2011JCP,kuchlin2018stochastic}, the unified stochastic particle method~\cite{fei2020usp,fei2021hybrid,fei2022usppoly}, and the general synthetic iterative scheme (GSIS)~\cite{SuArXiv2019,Su2020SIAM,Zhu2021JCP}.

These multiscale methods have also been extended to gas mixtures. The USP method has been developed for polyatomic gas mixtures~\cite{fei2026uspmix}. The discrete unified gas kinetic scheme has been applied to gas separation in pressure-driven binary-mixture flows through a microchannel~\cite{zhang2018discrete,zhang2019dugksmixture,xin2025dugkssep}. For the GSIS, a binary-mixture formulation~\cite{zeng2024gsismixture} has been developed based on the Li kinetic model~\cite{li2024kinetic}. However, the existing GSIS for gas mixtures employs a low-order finite-volume discretization. The asymptotic-preserving (AP) property removes the kinetic-scale restriction on the spatial cell size, but does not by itself reduce the spatial resolution required for accurately resolving macroscopic flow structures. Consequently, although the cell size is no longer constrained by the molecular mean free path in the near-continuum regime, a large number of physical cells may still be required for spatial accuracy when low-order discretizations are employed. For gas mixtures, this resolution requirement is more demanding because the flow structures of different species can have distinct spatial scales. Meanwhile, discretization in physical and molecular velocity spaces for multiple species leads to a large total number of degrees of freedom, increasing computational cost and memory requirements.

A natural approach to overcoming this remaining spatial-resolution bottleneck is to employ high-order spatial discretization, for which fewer physical elements are generally required to achieve a prescribed accuracy~\cite{su2015rkdg,su2020implicitdg}. Among various high-order methods, the discontinuous Galerkin (DG) method~\cite{cockburn2001rkdg} is particularly attractive for kinetic equations. High-order accuracy is obtained by increasing the polynomial degree within each element while retaining a compact stencil that couples face-neighboring elements~\cite{cockburn2001rkdg,su2015rkdg}. High-order discretizations have been extensively explored for kinetic equations, including WENO, spectral, and DG methods~\cite{filbet2003weno,gamba2009spectral,alekseenko2014dg}, while gas-kinetic DG and high-order compact schemes have also been developed for continuum flows~\cite{liu2007rkdg,ji2024compactgks}.High-order asymptotic-preserving DG schemes have also been developed for kinetic equations under different asymptotic scalings~\cite{xiong2015apdg,jang2015apdg}. For nonlinear kinetic equations, physical-space DG methods have been implemented with explicit Runge--Kutta time marching~\cite{su2015rkdg}, where the admissible time step becomes increasingly restrictive with the polynomial degree and is further limited by the molecular collision time. An implicit DG solver based on the conventional iteration scheme~\cite{su2020implicitdg} has demonstrated that high-order spatial discretization can substantially reduce the number of physical elements, but its convergence deteriorates in the near-continuum regime. On the other hand, synthetic iterations have recently been combined with high-order discretizations for several linear kinetic models~\cite{su2019hdg,liu2026gsiseph,shen2025agefem,shen2026nongray}, showing that fast synthetic convergence and high-order spatial accuracy can be achieved simultaneously. However, the combination of synthetic iteration and high-order DG discretization has not yet been developed for nonlinear kinetic equations of rarefied gas mixtures.

The aim of the present work is therefore to develop the GSIS-DG, namely, a high-order DG implementation of the GSIS, for efficient steady-state simulation of rarefied gas mixture flows over the entire range of Knudsen numbers. Based on the Li kinetic model~\cite{li2024kinetic} and the GSIS for gas mixtures~\cite{zeng2024gsismixture}, the mesoscopic kinetic equations and macroscopic synthetic equations are discretized by the DG method with the same polynomial order on the same physical mesh and solved alternately. The resulting method addresses two spatial-scale requirements in a complementary manner: the AP property of the GSIS removes the need to resolve the molecular mean free path in the near-continuum regime, while the high-order DG discretization reduces the number of physical elements required to resolve macroscopic flow structures accurately. Meanwhile, the synthetic equations retain the fast convergence of the GSIS toward the steady state. The accuracy, high-order spatial convergence, AP property, and computational efficiency of the GSIS-DG are assessed through normal shock waves, lid-driven cavity flows, and three-dimensional nozzle expansion, covering a wide range of gas rarefaction and flow configurations.

The remainder of this paper is organized as follows. Section~\ref{sec:2} introduces the Li kinetic model for gas mixtures and constructs the macroscopic synthetic equations in the GSIS. Section~\ref{sec:ap} analyzes continuous fixed-point compatibility and the formal Navier--Stokes limit, together with a sufficient spatial accuracy condition. Section~\ref{sec:3} then presents the DG discretizations of the kinetic and synthetic equations, the boundary treatment, and the overall GSIS-DG algorithm. The accuracy and efficiency of the GSIS-DG are assessed through three benchmark problems in Section~\ref{sec:4}. Finally, conclusions are drawn in Section~\ref{sec:5}.

\section{Kinetic model and macroscopic synthetic equations}\label{sec:2}

The governing equations consist of the Li kinetic model for a binary monatomic mixture~\cite{li2024kinetic} and the macroscopic synthetic equations used to accelerate its steady-state solution~\cite{Zhu2021JCP,zeng2024gsismixture}. We first specify the kinetic model and its parameters, and then define the reference constitutive relations and kinetic corrections used in the coupling.

\subsection{Li kinetic model}
 
The reference length, temperature, number density and molecular mass are denoted by $L_0$, $T_0$, $n_0$ and $m_0$, respectively. The corresponding velocity and time scales are $v_0=\sqrt{k_BT_0/m_0}$ and $L_0/v_0$, where $k_B$ is the Boltzmann constant. Using dimensional reference properties, the Knudsen number of species $s$ is defined as
\begin{equation}\label{eq:Kn_define}
\text{Kn}_{s} = \frac{\mu_s(T_0)}{n_0 L_0} \sqrt{\frac{\pi}{2m_s k_B T_0} },
\end{equation}
where $m_s$ is the molecular mass and $\mu_s(T_0)$ is the pure-species viscosity at the reference temperature. The dimensional viscosity follows $\mu_s(T_s)=\mu_s(T_0)(T_s/T_0)^{\omega_s}$, with $\omega_s=1$ for Maxwell molecules and $\omega_s=1/2$ for hard spheres. The mass density, pressure and heat flux are normalized by $m_0n_0$, $n_0k_BT_0$ and $n_0k_BT_0v_0$, respectively, and the mass distribution function by $m_0n_0/v_0^3$. All variables below are dimensionless, with the original symbols retained. Let $f_s(\bm x,\bm\xi,t)$ be the mass distribution function, whose velocity integral gives the mass density of species $s$. With the peculiar velocity $\bm c_s=\bm\xi-\bm u_s$, the species density, velocity, temperature, pressure tensor and heat flux are obtained from
\begin{equation}\label{eq:getmoment}
\begin{aligned}
    &\rho_s=\langle 1,f_s\rangle,\qquad n_s=\rho_s/m_s,\qquad
    \rho_s\bm u_s=\langle\bm\xi,f_s\rangle,\\
    &\frac32n_sT_s=\left\langle\frac12|\bm c_s|^2,f_s\right\rangle,\quad
    \bm P_s=\langle\bm c_s\bm c_s,f_s\rangle,\quad
    \bm q_s=\left\langle\frac12|\bm c_s|^2\bm c_s,f_s\right\rangle,
\end{aligned}
\end{equation}
where $\langle h,\psi\rangle=\int_{\mathbb R^3}h\psi\,\mathrm d\bm\xi$. Molecular velocities have three components irrespective of the physical-space dimension. The species pressure is $p_s=n_sT_s$, and the deviatoric stress is $\bm\sigma_s=\bm P_s-p_s\bm I$. The mixture number density, mass density and barycentric velocity satisfy
\begin{equation}
\begin{aligned}
n &= \sum_s n_s,\qquad \rho=\sum_s\rho_s,\qquad
\rho\bm u=\sum_s\rho_s\bm u_s.
\end{aligned}
\end{equation}
The mixture temperature, pressure tensor and heat flux are defined relative to $\bm u$. Introducing the diffusion velocity $\bm w_s=\bm u_s-\bm u$, with $\sum_s\rho_s\bm w_s=0$, their expressions in terms of the species moments are
\begin{equation}\label{eq:mixture-moments}
\begin{aligned}
\frac32nT&=\sum_s\frac32n_sT_s+\frac12\sum_s\rho_s|\bm w_s|^2,\\
\bm P&=\sum_s\bm P_s+\sum_s\rho_s\bm w_s\bm w_s,\\
\bm q&=\sum_s\bm q_s+\frac32\sum_s n_sT_s\bm w_s
+\frac12\sum_s\rho_s|\bm w_s|^2\bm w_s+\sum_s\bm P_s\cdot\bm w_s.
\end{aligned}
\end{equation}

For $s,r\in\{1,2\}$, the evolution of $f_s$ is governed by the Li kinetic equation~\cite{li2024kinetic}:
\begin{equation}\label{eq:generalModel}
\frac{\partial f_s}{\partial t}+\bm{\xi}\cdot \frac{\partial f_s}{\partial \bm{x}} = \sum_{r} \frac{g_{sr}-f_s}{\tau_{sr}}, 
\end{equation}
The terms with $r=s$ and $r\ne s$ describe intra- and inter-species collisions, respectively. Their relaxation times are related to the species Knudsen numbers through
\begin{equation}
     \tau_{ss}=\text{Kn}_{s}\sqrt{\frac{2 m_s}{\pi}}\frac{{T_s}^{\omega_s-1}}{{n}_s}, \qquad \tau_{sr} = \tau_{ss}\phi_{sr}^{-1}\frac{ n_s}{ n_r}.
\end{equation}
where parameter $\phi_{sr}$ controls the relative collision rates and is determined from the mixture viscosity. To write the reference distribution compactly, let $\hat{\bm c}_{sr}=\bm\xi-\hat{\bm u}_{sr}$. Then
\begin{equation}\label{reference distribution function}
\begin{aligned}
g_{sr}={}&m_s\hat n_{sr}\left(\frac{m_s}{2\pi T_s}\right)^{3/2}
\exp\left(-\frac{m_s|\hat{\bm c}_{sr}|^2}{2T_s}\right)\\
&\times\left[1+\frac{\hat T_{sr}-T_s}{T_s}
\left(\frac{m_s|\hat{\bm c}_{sr}|^2}{2T_s}-\frac32\right)
+\frac{2m_s\hat{\bm q}_{sr}\cdot\hat{\bm c}_{sr}}{5\hat n_{sr}T_s^2}
\left(\frac{m_s|\hat{\bm c}_{sr}|^2}{2T_s}-\frac52\right)\right].
\end{aligned}
\end{equation}
Here, $\hat n_{sr}$, $\hat{\bm u}_{sr}$, $\hat T_{sr}$ and $\hat{\bm q}_{sr}$ specify the reference moments. The Gaussian and both correction polynomials use $T_s$; the temperature correction sets the energy moment to $\hat T_{sr}$, while the heat-flux correction leaves the mass, momentum and energy moments unchanged. Conservation of species mass, momentum and energy in self-collisions requires
\begin{equation}\label{eq:aux_para_intra}
   \hat n_{ss} = n_s,\quad
    \hat {\bm{u}}_{ss} = {\bm{u}}_s,\quad
    \hat T_{ss} = T_s.
\end{equation}
For cross-collisions ($r\ne s$), the auxiliary states are
\begin{equation}\label{eq:aux_para_inter}
\begin{aligned}
\hat n_{sr} &= n_s, \\
\hat {\bm{u}}_{sr} &= \bu_s - \frac{\rho_r \tau_{sr}}{\rho_s \tau_{rs} + \rho_r\tau_{sr} }\bm{X}_{sr},\\
\hat {T}_{sr} &= T_s - \frac{n_r \tau_{sr}}{n_s \tau_{rs} + n_r\tau_{sr} }Y_{sr} - \frac{\rho_s \rho_r \tau_{sr} \tau_{rs} \bm{X}_{sr}\cdot [\bm{X}_{sr} - 2(\bm u_s - \bm u_r)]}{3(n_s\tau_{rs} + n_r \tau_{sr})(\rho_s \tau_{rs} + \rho_r\tau_{sr})},\\
\end{aligned}
\end{equation}
with
\begin{equation}
\begin{aligned}
\bm{X}_{sr} & = a_{sr}(\bm{u}_s - \bm{u}_r) + b_{sr} (\nabla \ln T_s + \nabla \ln T_r),\\
Y_{sr} & = c_{sr}(T_s - T_r) + d_{sr} |\bm{u}_s - \bm{u}_r|^2,\\
\end{aligned}
\end{equation}
The transport parameters of the Li model are~\cite{li2024kinetic}
\begin{equation}
    \begin{aligned}
        a_{sr} &= \frac{T\left(\rho_s\tau_{rs}+\rho_r\tau_{sr}\right)}{m_sm_r(n_s+n_r)D_{sr}}, \\
        b_{sr} &= \frac{(n_s+n_r)T\left(\rho_s\tau_{rs} + \rho_r\tau_{sr}\right)k_{T,sr}}{2\rho_s\rho_r}, \\
        c_{sr} &= \frac{2a_{sr}\left(n_s\tau_{rs}+n_r\tau_{sr}\right)m_sm_r}{(m_s+m_r)\left(\rho_s\tau_{rs}+\rho_r\tau_{sr}\right)}, \\
        d_{sr} &= \frac{a_{sr}m_sm_r}{3\left(\rho_s\tau_{rs}+\rho_r\tau_{sr}\right)}\left[\frac{a_{sr}\left(n_r\rho_r\tau_{sr}^2-n_s\rho_s\tau_{rs}^2\right)}{\rho_s\tau_{rs}+\rho_r\tau_{sr}} - \frac{2\left(\rho_r\tau_{sr}-\rho_s\tau_{rs}\right)}{m_s+m_r}\right],
    \end{aligned}
\end{equation}
where $D_{sr}$ and $k_{T,sr}$ are the binary diffusion coefficient and thermal-diffusion ratio. These parameters satisfy $a_{sr}=a_{rs}$, $b_{sr}=-b_{rs}$, $c_{sr}=c_{rs}$ and $d_{sr}=-d_{rs}$.

The auxiliary heat fluxes complete the reference distributions:
\begin{equation}\label{eq:aux_heatflux}
    \begin{aligned}
        \hat{\bm{q}}_{ss} &= (1-\text{Pr}) \bq_s, \\
        \hat{\bm{q}}_{sr} &= (1-\varphi_{sr}\text{Pr} ) \bq_s + \gamma_{sr} (\bq_{sr}-\bq_s), \\
        \bm q_{sr} &= \left\langle \frac{1}{2}|\bxi-\bm{\hat u}_{sr}|^2 (\bxi-\bm{\hat u}_{sr}), f_s \right\rangle.
    \end{aligned}
\end{equation}
Here, $\Pr$ is the Prandtl number, and $\varphi_{sr}$ adjusts the inter-species thermal relaxation rate to match the mixture thermal conductivity. The moment $\bm q_{sr}$ is evaluated relative to $\hat{\bm u}_{sr}$, whereas $\bm q_s$ is evaluated relative to $\bm u_s$. The Dufour coefficient $\gamma_{sr}=-\gamma_{rs}$ is specified by
\begin{equation}
\begin{aligned}
\gamma_{sr} &=\left(\frac{1}{m_s}\frac{\tau_{ss}\tau_{sr}}{\varphi_{sr}\tau_{ss}+\tau_{sr}}+\frac{1}{m_r}\frac{\tau_{rr}\tau_{rs}}{\varphi_{rs}\tau_{rr}+\tau_{rs}}\right)^{-1}\frac{4b_{sr}\Pr}{5a_{sr}T},
\end{aligned}
\end{equation}

Thus $\phi_{sr}$ and $\varphi_{sr}$ set the viscous and thermal relaxation rates, while $a_{sr}$, $b_{sr}$ and $\gamma_{sr}$ determine diffusion and its coupling to heat transport. The gas parameters used in the numerical tests are given in Section~\ref{sec:4}.

\subsection{Synthetic equations and constitutive coupling}

The macroscopic synthetic equations combine species mass conservation with momentum and energy balances that account for interspecies exchange. Together with constitutive relations that recover the Navier--Stokes limit and kinetic corrections that retain nonequilibrium effects, they provide a macroscopic system for accelerating the kinetic iteration. Taking the mass, momentum and total-energy moments of Eq.~\eqref{eq:generalModel} gives the governing equations for the species density $\rho_s$, velocity $\bm u_s$ and temperature $T_s$:
\begin{equation}\label{eq:macro_equation}
\begin{aligned}
\frac{\partial{\rho_s}}{\partial{t}} + \nabla\cdot\left(\rho_s\bm{u}_s\right)  &= 0, \\
    \frac{\partial}{\partial{t}}\left(\rho_s\bm{u}_s\right) + \nabla\cdot\left(\rho_s\bm{u}_s\bm{u}_s\right) + \nabla\cdot (n_s T_s \bm{\mathrm{I}}+ \bm{\sigma}_s) &= \bm{Q}_s^M, \\
    \frac{\partial}{\partial{t}}\left(E_s\right) + \nabla\cdot\left(E_s\bm{u}_s\right) + \nabla\cdot\left[(n_s T_s \bm{\mathrm{I}}+ \bm{\sigma}_s)\cdot\bm{u}_s+\bm{q}_s\right] &= {Q}_s^E,
\end{aligned}
\end{equation}
Here, $E_s=\frac32n_sT_s+\frac12\rho_s|\bm u_s|^2$ is the species total-energy density. The cross-collision sources are
\begin{equation}\label{momentum_exchange}
    \begin{aligned}
        \bm{Q}_s^M = \sum_{r\ne s} \frac{\rho_s(\hat{\bm{u}}_{sr}-\bm{u}_s)}{\tau_{sr}},\quad
        Q_s^E = \sum_{r\ne s} \frac{E_{sr}-E_s}{\tau_{sr}},
    \end{aligned}
\end{equation}
where $E_{sr}=\frac32n_s\hat T_{sr}+\frac12\rho_s|\hat{\bm u}_{sr}|^2$. The auxiliary states ensure $\sum_s\bm Q_s^M=0$ and $\sum_sQ_s^E=0$, so that collisions conserve mixture momentum and energy. The synthetic equations require the stress and heat flux supplied by the kinetic solution. To obtain the synthetic equations, these quantities are decomposed into reference constitutive relations and kinetic corrections~\cite{Zhu2021JCP}:
\begin{equation}\label{eq:full_constitutive}
    \begin{aligned}
        \bm{\sigma}_s &= \underbrace{-\mu_s^{\ast} \left[\nabla\bm{u}_s+(\nabla\bm{u}_s)^{\mathrm{T}}-\frac{2}{3}(\nabla\cdot\bm{u}_s)\bm I\right]}_{ \bm{\sigma}^{\mathrm{ref}}_s} + \mathrm{HoT}_{\bm{\sigma}_s},\\
        \bm{q}_{s} &=\underbrace{-\kappa_s^{\ast}\nabla T_s
        +\Lambda_{sr}(\bm u_s-\bm u_r)}_{\bm q_s^{\mathrm{ref}}}
        +\mathrm{HoT}_{\bm q_s},\qquad r\ne s,
    \end{aligned}
\end{equation}
The reference heat flux includes thermal conduction and the first-order contribution associated with species diffusion. Its coefficients follow from expressing the Chapman--Enskog heat flux in terms of the relative species velocity~\cite{li2024kinetic}. Introducing the symmetric interspecies momentum-exchange coefficient $\mathcal A_{sr}=\rho_s\rho_r/(\rho_s\tau_{rs}+\rho_r\tau_{sr})$, the constitutive coefficients are
\begin{equation}\label{eq:complete_heat_coefficients}
\begin{aligned}
 &\nu_s=\tau_{ss}^{-1}+\tau_{sr}^{-1},\quad
 t_{q,s}=\left(\frac{\Pr}{\tau_{ss}}+\frac{\varphi_{sr}\Pr}{\tau_{sr}}\right)^{-1},\quad
 \Lambda_{sr}=\frac{5T}{2m_s}t_{q,s}\gamma_{sr}a_{sr}\mathcal A_{sr},\\
 &\mu_s^0=\frac{n_sT_s}{\nu_s},\qquad
 \kappa_s^0=\frac{5n_sT}{2m_s}t_{q,s},\qquad
 \kappa_s^{\mathrm{eff}}=\kappa_s^0-\frac{2b_{sr}}{a_{sr}T}\Lambda_{sr}.
\end{aligned}
\end{equation}
Here, $T$ is the mixture temperature in Eq.~\eqref{eq:mixture-moments}. The relaxation time $t_{q,s}$ includes both self- and cross-collision heat-flux relaxation. The correction $-2b_{sr}\Lambda_{sr}/(a_{sr}T)$ in $\kappa_s^{\mathrm{eff}}$ accompanies the explicit relative-velocity term; these two contributions form a single first-order heat-flux representation. The reference transport coefficients are
\begin{equation}\label{eq:nominal_NSF_constitutive}
    \begin{aligned}
        \mu_s^{\ast} & = \beta_{s,\mu}\mu_s^0,\qquad
        \beta_{s,\mu} = \frac{1}{\text{Kn}_s}\text{min}\left(\text{Kn}_s,\frac{1}{n}\sum_r n_r\phi_{sr}\right), \\
        \kappa_s^{\ast} &= \beta_{s,\kappa}\kappa_s^{\mathrm{eff}},\qquad
        \beta_{s,\kappa} = \frac{1}{\text{Kn}_s}\text{min}\left(\text{Kn}_s,\frac{1}{n}\sum_r n_r\phi_{sr}\varphi_{sr}\right),
    \end{aligned}
\end{equation}
with $\phi_{ss}=\varphi_{ss}=1$. The factors $\beta_{s,\mu}$ and $\beta_{s,\kappa}$ limit the gradient coefficients at large Knudsen numbers and become unity in the continuum limit, while $\Lambda_{sr}$ is not rescaled. The diffusion term and the correction in $\kappa_s^{\mathrm{eff}}$ together recover the first-order Li heat flux. Let $\bm W_s^f=(\rho_s^f,\rho_s^f(\bm u_s^f)^{\mathrm T},E_s^f)^{\mathrm T}$ denote the conservative moments of species $s$, and let $\mathbb W^f=(\bm W_1^f,\bm W_2^f)$ collect both species states. The higher-order terms (HoTs) are evaluated by subtracting the complete reference constitutive relations, including diffusion, at this kinetic moment state:
\begin{equation}\label{eq:getHoTs}
        \mathrm{HoT}_{\bm{\sigma}_s} = \bm\sigma_s[f_s]-\bm\sigma_s^{\mathrm{ref}}(\mathbb W^f),\quad
        \mathrm{HoT}_{\bm{q}_s} = \bm q_s[f_s]-\bm q_s^{\mathrm{ref}}(\mathbb W^f).
\end{equation}

Here, $\bm{\sigma}_s[f_s]=\bm{P}_s[f_s]-n_sT_s\bm{I}$ and $\bm{q}_s[f_s]$ are evaluated from the kinetic solution, while $\mathbb{W}^f$ is used to evaluate the corresponding reference constitutive terms. During the macroscopic inner iterations, the HoTs remain fixed, whereas the reference constitutive terms are updated with the synthetic state. When the synthetic state coincides with the kinetic moments, the reference terms cancel those subtracted in Eq.~\eqref{eq:getHoTs}, so that the kinetic stress and heat flux are exactly recovered. The fixed-point consistency and continuum limit are analyzed in Section~\ref{sec:ap}, and the DG discretization is described in Section~\ref{sec:3}.
\section{Asymptotic-preserving analysis}
\label{sec:ap}

This section examines the asymptotic behavior~\cite{jin2022ap} of the converged GSIS. Following the analysis of Shen \emph{et al.}~\cite{shen2026nongray}, we first establish the fixed-point consistency between the kinetic and synthetic equations, and then consider the Navier--Stokes (NS) limit and the corresponding spatial-accuracy requirement. The analysis assumes smooth near-equilibrium solutions with positive densities and temperatures, fixed molecular-mass and concentration ratios, and nondegenerate relaxation rates. Molecular-velocity moments are assumed to be evaluated exactly, and periodic or compatible boundary conditions are considered so that unresolved Knudsen layers do not enter the analysis.

\subsection{Fixed-point consistency}
\label{sec:ap-fixed-point}

At a fixed point of the GSIS, let $\bm{W}_s^*$ denote the conservative state obtained from the macroscopic synthetic equations, and let $\bm{W}_s^f$ denote the corresponding conservative moments of $f_s$. Their difference is written as
\begin{equation}
\delta\bm{W}_s
=
\bm{W}_s^*-\bm{W}_s^f
=
\left(
\delta\rho_s,
\delta(\rho_s\bm{u}_s)^{\mathrm T},
\delta E_s
\right)^{\mathrm T}.
\label{eq:ap-state}
\end{equation}
The two states are not assumed \emph{a priori} to coincide. At the fixed point, the low-order parameters of the collision reference distributions in Eq.~\eqref{eq:generalModel} are evaluated from the collective synthetic state $\mathbb W^*=(\bm W_1^*,\bm W_2^*)$. Collect the exchange terms in Eq.~\eqref{momentum_exchange} into the conservative source $\bm{Q}_s$, and let $\bm{F}_s^f$ and $\bm{F}_s^{\mathrm{syn}}$ denote the conservative fluxes of the kinetic moment and synthetic equations, respectively. The latter is evaluated using the kinetic HoTs defined in Eq.~\eqref{eq:getHoTs}. Taking conservative moments of the steady form of Eq.~\eqref{eq:generalModel} and comparing them with the synthetic equations gives
\begin{equation}
\nabla\cdot\bm{F}_s^f
=
\bm{Q}_s(\mathbb W^*)
+\nu_s\delta\bm{W}_s,
\quad
\nabla\cdot\bm{F}_s^{\mathrm{syn}}
=
\bm{Q}_s(\mathbb W^*),
\label{eq:ap-moment-equations}
\end{equation}
The additional term in the first equation arises because the collision-loss moment is evaluated from $\bm{W}_s^f$, whereas the collision reference distributions are constructed from $\bm{W}_s^*$. The interspecies momentum and energy exchanges are identical in the two equations and therefore cancel exactly when they are subtracted. Consequently, the difference between the synthetic state and the kinetic moments satisfies
\begin{equation}
\nu_s\delta\bm{W}_s
+
\nabla\cdot
\left[
\bm{F}_s^{\mathrm{syn}}
-
\bm{F}_s^f
\right]
=0.
\label{eq:ap-increment}
\end{equation}
Assuming that this homogeneous increment problem, together with the corresponding homogeneous boundary conditions, admits only the trivial solution, one has
$\delta\bm{W}_s=0$ and hence $\bm{W}_s^*=\bm{W}_s^f$.
At this fixed point, the reference stress and heat-flux contributions evaluated from the synthetic state are identical to those subtracted from the kinetic moments in the HoTs and therefore cancel exactly. Meanwhile, $g_{sr}(\mathbb W^*)$ recovers the collision reference distribution of the original Li model. Thus, the converged coupled system is consistent with the original steady kinetic equations; conversely, every steady solution of the original kinetic equations together with its moments satisfies the coupled continuous equations.

\subsection{Continuum limit and spatial accuracy}
\label{sec:ap-ns}

Let $\varepsilon$ denote a reference Knudsen number, with $\Kn_s=O(\varepsilon)$ and fixed ratios between the species Knudsen numbers. Under the continuum scaling of the Li model, $\tau_{sr}=O(\varepsilon)$ and $b_{sr}=O(\varepsilon)$, while the remaining nondimensional collision parameters are $O(1)$. For sufficiently small $\varepsilon$, the rescaling factors in Eq.~\eqref{eq:nominal_NSF_constitutive} satisfy
$\beta_{s,\mu}=\beta_{s,\kappa}=1$.

The point at which the mixture structure enters is the stiff interspecies momentum exchange. For $r\ne s$, let $\bm w_{sr}=\bm u_s-\bm u_r$. Substituting the auxiliary velocity in Eq.~\eqref{eq:aux_para_inter} into Eq.~\eqref{momentum_exchange} gives $\bm Q_s^M=-\mathcal A_{sr}\bm X_{sr}$. The coefficient $\mathcal A_{sr}$ sets the strength of interspecies momentum exchange and scales inversely with the collision relaxation time. Expanding $\bm X_{sr}$ yields
\begin{equation}
\bm Q_s^M
=-
\mathcal A_{sr}
\left[
a_{sr}\bm w_{sr}
+b_{sr}\left(\nabla\ln T_s+\nabla\ln T_r\right)
\right],
\qquad
\mathcal A_{sr}=O(\varepsilon^{-1}),
\quad
\Lambda_{sr}=O(1).
\label{eq:ap-stiff-momentum-exchange}
\end{equation}
An $O(1)$ species pressure-gradient imbalance is balanced by the $O(\varepsilon^{-1})$ exchange operator, giving $\bm w_{sr}=O(\varepsilon)$, with ordinary and thermal diffusion both entering at first order. The Chapman--Enskog expansion of the Li model~\cite{li2024kinetic} gives $T_s-T_r=O(\varepsilon^2)$ under the present scaling. From Eq.~\eqref{eq:complete_heat_coefficients}, $\mu_s^0$ and $\kappa_s^{\mathrm{eff}}$ are $O(\varepsilon)$, whereas $\Lambda_{sr}=O(1)$, so the gradient and relative-velocity contributions to the reference heat flux are both of first order. With $\beta_{s,\mu}=\beta_{s,\kappa}=1$, the reference relations in Eq.~\eqref{eq:full_constitutive} reproduce the complete first-order species stress and heat flux. Their subtraction in Eq.~\eqref{eq:getHoTs} therefore yields
\begin{equation}
\mathrm{HoT}_{\bm\sigma_s}
=
O(\varepsilon^2),
\qquad
\mathrm{HoT}_{\bm q_s}
=
O(\varepsilon^2).
\label{eq:ap-hot}
\end{equation}
Using Eq.~\eqref{eq:ap-hot} in the synthetic equations and summing the momentum and energy equations over the two species yields, through first order in $\varepsilon$,
\begin{align}
\nabla\cdot
\left(
\rho_s\bm{u}
+
\bm{J}_s^{\mathrm{CE}}
\right)
&=
O(\varepsilon^2),
\qquad s=1,2,
\nonumber\\
\nabla\cdot
\left(
\rho\bm{u}\bm{u}
+
p\bm{I}
+
\bm{\sigma}^{\mathrm{NS}}
\right)
&=
O(\varepsilon^2),
\nonumber\\
\nabla\cdot
\left[
(E+p)\bm{u}
+
\bm{\sigma}^{\mathrm{NS}}\cdot\bm{u}
+
\bm{q}^{\mathrm{NS}}
\right]
&=
O(\varepsilon^2),
\label{eq:ap-ns-limit}
\end{align}
where $p=nT$, $E=\sum_s E_s$, and $\bm J_s^{\mathrm{CE}}=\rho_s\bm w_s^{\mathrm{CE}}$ is the first-order diffusion mass flux, with $\sum_s\bm J_s^{\mathrm{CE}}=0$, while
$\bm{\sigma}^{\mathrm{NS}}$ and $\bm{q}^{\mathrm{NS}}$ are the corresponding mixture stress and heat flux derived in Ref.~\cite{li2024kinetic}. In particular, Eq.~\eqref{eq:mixture-moments} gives $\bm q^{\mathrm{NS}}=\sum_s\bm q_s^{\mathrm{ref}}+\frac52\sum_s n_sT\bm w_s$ through first order. The diffusion enthalpy arises from summing the species energy fluxes and is not counted again in $\bm q_s^{\mathrm{ref}}$. The interspecies momentum and energy exchange terms cancel by conservation when the species equations are summed. Hence, once fixed-point consistency is reached, the continuous synthetic equations recover the multicomponent NS equations in the limit $\varepsilon\rightarrow0$. Although $T_s-T_r=O(\varepsilon^2)$, the individual interspecies energy-exchange term need not be negligible: division by $\tau_{sr}=O(\varepsilon)$ can leave an $O(\varepsilon)$ contribution. It disappears from the mixture energy equation through pairwise conservation after summing the species equations, rather than by being discarded in each species equation.

It remains to connect the continuous limit to the spatial discretization. Let $h$ denote the maximum physical element diameter divided by $L_0$. At finite $h$, the synthetic state and the kinetic moments need not coincide exactly because their numerical residual operators are different. Assume that the complete kinetic and macroscopic DG operators are $K$th-order consistent for the smooth solution family considered here and that the normalized collision linearization is uniformly invertible on the nonconserved modes. A formal perturbation argument then gives a fixed-point state difference of $O(\varepsilon h^K)$. More specifically, if $\bm w_{sr,h}^{\mathrm{CE}}$ denotes the Li first-order relative velocity evaluated with consistent discrete gradients, then the $O(h^K)$ consistency defect in the species momentum balance is filtered through the factor $\mathcal A_{sr}^{-1}=O(\varepsilon)$ in Eq.~\eqref{eq:ap-stiff-momentum-exchange}. Consequently,
\begin{equation}
\bm w_{sr,h}-\bm w_{sr,h}^{\mathrm{CE}}
=O(\varepsilon^2)+O(\varepsilon h^K),
\qquad
\Lambda_{sr}
\left(\bm w_{sr,h}-\bm w_{sr,h}^{\mathrm{CE}}\right)
=O(\varepsilon^2)+O(\varepsilon h^K).
\label{eq:ap-relative-velocity-defect}
\end{equation}
The second estimate is essential for the diffusion-induced heat flux: it shows that this term does not introduce an unscaled $O(h^K)$ constitutive defect. Propagating Eq.~\eqref{eq:ap-relative-velocity-defect}, the remaining first-order constitutive defects, and one spatial derivative through the summed equations gives the conditional residual estimate
\begin{equation}
\left\|
\mathcal{R}_{\mathrm{NS,mix},h}^{\varepsilon}
\right\|
\leq
C\left(
\varepsilon^2+h^K+\varepsilon h^{K-1}
\right),
\label{eq:ap-discrete-residual}
\end{equation}
where $\mathcal{R}_{\mathrm{NS,mix},h}^{\varepsilon}$ is the residual of the complete first-order mixture NS system evaluated on a smooth reconstruction of the discrete state. This is a formal residual estimate rather than an error estimate relative to the exact NS solution. A sufficient condition for the two spatial terms in Eq.~\eqref{eq:ap-discrete-residual} to be no larger than the second-order Chapman--Enskog remainder, while leaving the molecular mean free path unresolved, is
\begin{equation}
h^K=O(\varepsilon^2),
\qquad
\frac{h}{\varepsilon}\rightarrow\infty.
\label{eq:ap-spatial-condition}
\end{equation}
Indeed, the first condition controls $h^K$, while the second gives $\varepsilon h^{K-1}=(\varepsilon/h)h^K=o(\varepsilon^2)$ and allows the physical mesh to remain much coarser than the kinetic length scale.

Writing $h\asymp\varepsilon^{1/\alpha}$, the two requirements become
\begin{equation}
h\asymp\varepsilon^{1/\alpha}
\quad\Longrightarrow\quad
1<\alpha\leq\frac{K}{2}
=\frac{p+1}{2},
\qquad K=p+1.
\label{eq:ap-dg-condition}
\end{equation}
Therefore, for sufficiently high polynomial order, the spatial truncation error can remain of second order or smaller in $\varepsilon$ even when the element size is asymptotically much larger than the molecular mean free path. Equation~\eqref{eq:ap-dg-condition} provides a sufficient rather than necessary condition and should not be interpreted as implying that lower-order discretizations are not asymptotic preserving.

\section{Numerical schemes of the GSIS-DG}\label{sec:3}

Following the nonlinear GSIS~\cite{Zhu2021JCP}, the kinetic and macroscopic synthetic equations are discretized on the same physical mesh and polynomial space. The superscripts $k$ and $k_{\mathrm{NS}}$ denote the outer GSIS and inner macroscopic iterations, respectively. The two solvers retain their own states between outer iterations. Inner states carry both indices $(k,k_{\mathrm{NS}})$, with $k$ held fixed during a macroscopic solve. Species are indexed by $s,r$, elements by $i,j$, and local basis functions by $a,b$. Their spatial discretizations, boundary treatments, and coupling procedure are presented below.

\subsection{DG discretization of the kinetic equation}\label{sec:dg-cis-scheme}

For notational simplicity, we write $f$ for $f_s$, which applies to both species. The collision operator is written as $\nu(g-f)$, with $\nu=\nu_s$ as defined in Eq.~\eqref{eq:complete_heat_coefficients} and $g=\nu^{-1}\sum_{s'=1}^2g_{ss'}/\tau_{ss'}$. Here, $s,s'\in\{1,2\}$ label species, so both self- and interspecies collisions are retained. At molecular velocity $\bm\xi_m$, the kinetic equation reads $\partial_t f_m+\bm\xi_m\cdot\nabla f_m=\nu(g_m-f_m)$, with $t$ denoting the pseudo time used to obtain the steady solution. The physical domain is partitioned into non-overlapping elements $\Omega_i$. Following the DG formulation of Su et al.~\cite{su2020implicitdg}, the distribution function within each element is approximated using Lagrange polynomials of degree $p$:
\begin{equation}\label{eq:dg-vdf-expansion}
 f_{h,m}\big|_{\Omega_i}
 =\sum_{b=1}^{N_p}f_{ibm}\ell_{ib}.
\end{equation}
where $\ell_{ib}$ are the local basis functions and $f_{ibm}$ are the corresponding expansion coefficients. For tensor-product elements in $d$ spatial dimensions, the number of coefficients is $N_p=(p+1)^d$. Multiplying the kinetic equation by the test function $\ell_{ia}$ and integrating the streaming term by parts gives the weak DG form, for $a=1,\ldots,N_p$:
\begingroup
\small
\begin{equation}\label{eq:dg-kinetic-weak}
 \sum_{b=1}^{N_p}(\bm M_i)_{ab}\frac{\mathrm d f_{ibm}}{\mathrm dt}
 =\int_{\Omega_i}(\bm\xi_m\cdot\nabla\ell_{ia})f_{h,m}\,\mathrm d\bm x
 +\int_{\Omega_i}\ell_{ia}\nu(g_m-f_{h,m})\,\mathrm d\bm x
 -\sum_{\Gamma\subset\partial\Omega_i}
 \int_\Gamma\ell_{ia}\widehat F_m\,\mathrm dS.
\end{equation}
\endgroup
Here, the mass matrix has entries $(\bm M_i)_{ab}=\int_{\Omega_i}\ell_{ia}\ell_{ib}\,\mathrm d\bm x$. The three right-hand-side terms represent volume transport, collisions, and outward convective flux, respectively, and together define the mesoscopic residual component $R^f_{iam}$. Since the VDF is discontinuous between elements, the numerical flux depends on the interior and exterior traces, denoted by $-$ and $+$. With $\bm n$ pointing out of the current element, the upwind flux reads
\begin{equation}\label{eq:dg-kinetic-upwind}
 \widehat F_m=\xi_n^+f_{h,m}^-+\xi_n^-f_{h,m}^+,
 \qquad \xi_n=\bm\xi_m\cdot\bm n,\qquad
 \xi_n^\pm=\frac{\xi_n\pm|\xi_n|}{2}.
\end{equation}
On a physical boundary, the incoming trace is prescribed by the boundary condition in Section~\ref{sec:dg-boundary-scheme}. At iteration $k$, the collision frequency and the low-order auxiliary parameters of the reference distributions are evaluated from the retained synthetic state $\mathbb W_h^{*,k}$. The nonequilibrium heat-flux inputs are obtained from $f_h^k$ using the moment definitions in Section~\ref{sec:2}. These quantities are held fixed during the kinetic update. When the synthetic state agrees with the conservative moments of $f_h$, this construction recovers the original Li collision operator. Backward Euler discretization in pseudo time treats transport and collision loss implicitly. For each element and molecular velocity, we collect the expansion coefficients into $\bm f_{i,m}=(f_{i1m},\ldots,f_{i{N_p}m})^{\mathrm T}$ and the residual components into $\bm R^f_{i,m}$. Introducing $\Delta\bm f_{i,m}=\bm f_{i,m}^{k+1}-\bm f_{i,m}^k$, the implicit update takes the element-block form
\begin{equation}\label{eq:delta_form_micro}
 \bm A_{ii,m}\Delta\bm f_{i,m}
 +\sum_{j\in N(i)}\bm A_{ij,m}\Delta\bm f_{j,m}
 =\bm R^{f,k}_{i,m}.
\end{equation}
where $N(i)$ denotes the face neighbors of element $i$ and $\Delta t_i$ is the local pseudo time step. Each block has size $N_p\times N_p$. With incoming physical-boundary traces held fixed during this update, the blocks are
\begin{equation}\label{eq:dg-kinetic-blocks}
\begin{aligned}
 \bm A_{ij,m}
 &=\frac{\delta_{ij}\bm M_i}{\Delta t_i}
 -\frac{\partial\bm R^f_{i,m}}{\partial\bm f_{j,m}},\qquad
 \bm A_{ii,m}=\frac{\bm M_i}{\Delta t_i}
 +\bm C_i-\bm J^{\mathrm{vol}}_{i,m}+\bm J^{\mathrm{face}}_{i,m},\\
  (\bm A_{ij,m})_{ab}
 &=\int_{\Gamma_{ij}}\xi_n^-\ell_{ia}\ell_{jb}\,\mathrm dS,\qquad
 (\bm C_i)_{ab}=\int_{\Omega_i}\nu\ell_{ia}\ell_{ib}\,\mathrm d\bm x,\\
 (\bm J^{\mathrm{vol}}_{i,m})_{ab}
 &=\int_{\Omega_i}(\bm\xi_m\cdot\nabla\ell_{ia})\ell_{ib}\,\mathrm d\bm x,\qquad
 (\bm J^{\mathrm{face}}_{i,m})_{ab}
 =\sum_{\Gamma\subset\partial\Omega_i}
 \int_\Gamma\xi_n^+\ell_{ia}\ell_{ib}\,\mathrm dS.
\end{aligned}
\end{equation}

The off-diagonal expression applies to $j\in N(i)$, $j\ne i$, with $\Gamma_{ij}$ denoting the shared face. For each molecular velocity, the elements are ordered along its upwind transport direction, and the block system is approximately solved using forward and backward LU-SGS sweeps:
\begin{equation}\label{eq:dg-kinetic-sweeps}
\begin{aligned}
 \Delta\bm f_{i,m}^{\mathrm{fwd}}
 &=\bm A_{ii,m}^{-1}
 \left(\bm R^{f,k}_{i,m}-\sum_{j<i}\bm A_{ij,m}\Delta\bm f_{j,m}^{\mathrm{fwd}}\right),\\
 \Delta\bm f_{i,m}
 &=\Delta\bm f_{i,m}^{\mathrm{fwd}}
 -\bm A_{ii,m}^{-1}\sum_{j>i}\bm A_{ij,m}\Delta\bm f_{j,m}.
\end{aligned}
\end{equation}
The sums run over face neighbors with $j<i$ and $j>i$ in this velocity-dependent ordering~\cite{su2020implicitdg}, respectively. The resulting increment gives $\bm f_{i,m}^{k+1}=\bm f_{i,m}^k+\Delta\bm f_{i,m}$. Its moments are then used to construct the kinetic constitutive corrections, without resetting the retained synthetic state.

\subsection{DG discretization of the macroscopic synthetic equations}\label{sec:dg-macro-scheme}

During the macroscopic inner iterations, the kinetic constitutive corrections remain fixed, while the reference constitutive terms and interspecies exchange terms are updated. We suppress the species index in the local expansion and residual only; $\bm W$ then denotes a single-species state, and both species are solved together with their coupling retained. With $\bm W^*=(\rho,\rho\bm u^{\mathrm T},E)^{\mathrm T}$, the local approximation is
\begin{equation}\label{eq:dg-projection}
 \bm W_h^{*,k,k_{\mathrm{NS}}}\big|_{\Omega_i}
 =\sum_{b=1}^{N_p}\ell_{ib}\bm W_{ib}^{*,k,k_{\mathrm{NS}}}.
\end{equation}
Here, $\bm W_{ib}^*$ contains the conservative-variable coefficients associated with basis function $\ell_{ib}$. The initial inner state is the previous synthetic solution, $\mathbb W_h^{*,k,0}=\mathbb W_h^{*,k}$, rather than the latest kinetic moments. Let $\bm W_{s,h}^{f,k+1}=\Pi_h\langle\bm\psi_d,f_{s,h}^{k+1}\rangle$, where $\Pi_h$ is the elementwise mass projection and $\bm\psi_d=(1,\xi_1,\ldots,\xi_d,|\bm\xi|^2/2)^{\mathrm T}$. For the planar and one-dimensional flows considered here, the remaining mean-velocity components vanish by symmetry; only the $d$ physical-space momentum components are retained, but the energy moment includes all three molecular-velocity components. Thus each species state has $d+2$ components. Velocity reduction, when used, must retain the transverse-energy moment. We write $\mathbb W_h=(\bm W_{1,h},\bm W_{2,h})$ for the collective state required by the coupled constitutive terms. The frozen constitutive corrections are evaluated at this kinetic moment state:
\begin{equation}\label{eq:dg-hot-projection}
\begin{aligned}
 \mathrm{HoT}_{\bm\sigma,h}^{k+1}
 &=\Pi_h\left[\bm\sigma^{k+1}
 -\bm\sigma^{\mathrm{ref}}(\mathbb W_h^{f,k+1})\right],\\
 \mathrm{HoT}_{\bm q,h}^{k+1}
 &=\Pi_h\left[\bm q^{k+1}
 -\bm q^{\mathrm{ref}}(\mathbb W_h^{f,k+1})\right].
\end{aligned}
\end{equation}
The projection satisfies $\sum_{b=1}^{N_p}(\bm M_i)_{ab}(\Pi_hz)_{ib}=\sum_q\omega_{i,q}\ell_{ia}(\bm x_{i,q})z(\bm x_{i,q})$, componentwise for vector and tensor fields. Here, $q$ indexes the volume quadrature points, and $\omega_{i,q}$ includes the physical volume Jacobian. The complete reference constitutive operator in Eqs.~\eqref{eq:full_constitutive}--\eqref{eq:nominal_NSF_constitutive}, including the diffusion heat flux, uses the reconstructed gradients of the synthetic solver and is evaluated consistently in the kinetic-state subtraction and the synthetic-state flux. The total flux is $\bm F=\bm F_c+\bm F_v$, with the convective and constitutive contributions given by Eqs.~\eqref{eq:macro_equation} and~\eqref{eq:full_constitutive}. During the inner iterations, the projected corrections remain fixed and are evaluated at the volume quadrature points. They are added to the current reference stresses and heat fluxes to assemble the volume flux in Eq.~\eqref{eq:dg-macro-residual}. Testing the synthetic equations with $\ell_{ia}$ and integrating the flux divergence by parts gives the weak DG form, for $a=1,\ldots,N_p$:
\begingroup
\small
\begin{equation}\label{eq:dg-macro-residual}
 \sum_{b=1}^{N_p}(\bm M_i)_{ab}\frac{\mathrm d\bm W_{ib}^*}{\mathrm dt}
 =\int_{\Omega_i}\nabla\ell_{ia}\cdot\bm F\,\mathrm d\bm x
 -\sum_{\Gamma\subset\partial\Omega_i}
 \int_\Gamma\ell_{ia}\widehat{\bm F}_n\,\mathrm dS
 +\int_{\Omega_i}\ell_{ia}\bm Q\,\mathrm d\bm x.
\end{equation}
\endgroup
The three right-hand-side terms together define the macroscopic residual $\bm R^{\mathrm{NS}}_{ia}$. Here, $\bm Q=(0,(\bm Q^M)^{\mathrm T},Q^E)^{\mathrm T}$ contains the interspecies momentum and energy exchange terms defined in Section~\ref{sec:2}. The convective and gradient-dependent viscous numerical fluxes are evaluated using the AUSM+up and BR2 schemes, respectively; the algebraic relative-velocity contribution to the reference heat flux is evaluated from the corresponding macroscopic states. The kinetic stresses and heat fluxes are first projected onto the DG space within each element. Their traces from adjacent elements are then arithmetically averaged at interior faces. Subtracting the complete reference face contributions evaluated at $\mathbb W_h^{f,k+1}$ gives the face corrections, which remain fixed during the inner iterations. The same face operators are used for this subtraction and for adding the reference contributions at the current synthetic state. This matches the constitutive operators, but does not identify the macroscopic numerical flux with the velocity moment of the upwind kinetic flux. For the implicit update, we collect the local coefficients into $\bm W_i^*=((\bm W_{i1}^*)^{\mathrm T},\ldots,(\bm W_{i{N_p}}^*)^{\mathrm T})^{\mathrm T}$ and the residuals into $\bm R_i^{\mathrm{NS}}$ in the same order. The corresponding block mass matrix is $\bm M_i^W=\bm M_i\otimes\bm I_{d+2}$. Restoring the species indices, we define $\Delta\bm W_{s,i}^*=\bm W_{s,i}^{*,k,k_{\mathrm{NS}}+1}-\bm W_{s,i}^{*,k,k_{\mathrm{NS}}}$. The kinetic corrections remain fixed, while the exchange sources are evaluated from the current states of both species. With local pseudo time step $\Delta t_i^{\mathrm{NS}}$, the coupled increment system is written as
\begin{equation}\label{eq:dis_macro}
 \sum_{s'=1}^2\sum_j\bm B_{ij}^{ss'}\Delta\bm W_{s',j}^*
 =\bm R_{s,i}^{\mathrm{NS},k,k_{\mathrm{NS}}},\qquad
 \bm B_{ij}^{ss'}=
 \frac{\delta_{ss'}\delta_{ij}\bm M_i^W}{\Delta t_i^{\mathrm{NS}}}
 -\frac{\partial\bm R_{s,i}^{\mathrm{NS}}}{\partial\bm W_{s',j}^*}.
\end{equation}
The Jacobian is evaluated at the current inner state with the kinetic corrections fixed, retaining spatial transport and interspecies source derivatives. The constitutive linearization includes the state dependence of the reference transport coefficients and the dependence of the diffusion heat flux on both species velocities. The index $j$ includes elements coupled through the gradient lifting. GMRES solves the system with a cell-block preconditioner that retains the coupling between both species. The state is then updated as $\bm W_{s,i}^{*,k,k_{\mathrm{NS}}+1}=\bm W_{s,i}^{*,k,k_{\mathrm{NS}}}+\Delta\bm W_{s,i}^*$.

\subsection{Boundary treatment}\label{sec:dg-boundary-scheme}

The physical boundary conditions are imposed through the incoming distribution at each face quadrature point; the outgoing distribution is reconstructed from the interior DG coefficients. We restore the species index $s$ and retain $m$ for the molecular-velocity index. At a diffuse wall, $\bm n$ points out of the gas domain and $\xi_n=(\bm\xi_{s,m}-\bm u_w)\cdot\bm n$, with $\bm u_w\cdot\bm n=0$ for the walls considered here. Let $M_s(\bm\xi;\bm u,T)=(m_s/(2\pi T))^{3/2}\exp[-m_s|\bm\xi-\bm u|^2/(2T)]$ denote the unit-density Maxwellian. Zero normal mass flux determines the reflected density separately at every face quadrature point:
\begin{equation}\label{eq:dg-wall-vdf}
\begin{aligned}
 f_{s,w}(\bm\xi_{s,m})&=
 \begin{cases}
 f_{s,h}^-(\bm\xi_{s,m}),&\xi_n\geq0,\\
 \rho_{s,w}M_s(\bm\xi_{s,m};\bm u_w,T_w),&\xi_n<0,
 \end{cases}\\
 \rho_{s,w}&=
 -\frac{\sum_{\xi_n>0}w_{s,m}\xi_n f_{s,h}^-(\bm\xi_{s,m})}
 {\sum_{\xi_n<0}w_{s,m}\xi_n M_s(\bm\xi_{s,m};\bm u_w,T_w)}.
\end{aligned}
\end{equation}
Here, $w_{s,m}$ denotes the molecular-velocity quadrature weight associated with $\bm\xi_{s,m}$. At a reservoir, the incoming distribution is the prescribed equilibrium distribution; at a vacuum outlet, it is zero. The boundary parameters are specified in the individual numerical examples.

The synthetic boundary flux is corrected by the moments of the latest kinetic boundary distribution~\cite{Zhu2021JCP,zeng2024gsismixture}. Let $\widehat{\bm F}_{s,b}^{M}(\mathbb W)$ denote the flux formed from outgoing and incoming Maxwellians satisfying the physical boundary condition. The difference between the kinetic flux and this Maxwellian flux is evaluated at $\mathbb W_h^{f,k+1}$ and fixed during the inner iterations. Thus, the corrected boundary flux agrees with the kinetic flux when evaluated at the kinetic moment state, while its Maxwellian part responds to the evolving synthetic state. Using the moment vector $\bm\psi_d$ defined in Section~\ref{sec:dg-macro-scheme}, the construction at each face quadrature point is
\begin{equation}\label{eq:dg-boundary-correction}
 \widehat{\bm F}_{s,b}^{k,k_{\mathrm{NS}}}
 =\widehat{\bm F}_{s,b}^{M}(\mathbb W_h^{*,k,k_{\mathrm{NS}}})
 +[\widehat{\bm F}_{s,b}^{\mathrm{kin}}
 -\widehat{\bm F}_{s,b}^{M}(\mathbb W_h^{f,k+1})], \quad   \widehat{\bm F}_{s,b}^{\mathrm{kin}}
 =\sum_m w_{s,m}\bm\psi_d(\bm\xi_{s,m})\xi_n
 f_{s,b}^{k+1}(\bm\xi_{s,m}).
\end{equation}
This total flux supplies the boundary integral in Eq.~\eqref{eq:dg-macro-residual}, without adding a separate viscous boundary flux.

\subsection{Overview of the GSIS-DG algorithm}\label{sec:dg-algorithm}

At outer iteration $k$, suppose that the kinetic state $f_h^k$ and the synthetic state $\mathbb W_h^{*,k}$ are available. For notational simplicity, the species index $s$ is omitted in this subsection; the same procedure is applied to both species with their interspecies coupling retained.
\begin{enumerate}
 \item Construct the collision reference distributions $g_{sr}(\mathbb W_h^{*,k})$. Solve Eq.~\eqref{eq:delta_form_micro} to advance the distribution function from $f_h^k$ to $f_h^{k+1}$.
 \item Use $f_h^{k+1}$ to evaluate the kinetic HoTs in Eq.~\eqref{eq:dg-hot-projection} and the boundary-flux correction in Eq.~\eqref{eq:dg-boundary-correction}. Starting from $\mathbb W_h^{*,k,0}=\mathbb W_h^{*,k}$, solve Eq.~\eqref{eq:dg-macro-residual} iteratively until the prescribed maximum number of inner iterations is reached or the macroscopic residual falls below the specified tolerance. If $L_k$ inner updates are performed, set $\mathbb W_h^{*,k+1}=\mathbb W_h^{*,k,L_k}$.
\end{enumerate}
The cycle is repeated with $(f_h^{k+1},\mathbb W_h^{*,k+1})$ until convergence. The updated synthetic state enters the collision reference distributions in the next kinetic solve; no equilibrium-distribution difference is added directly to the VDF. Evaluating the collision parameters from the moments of $f_h$ and omitting the macroscopic step recovers CIS-DG, whereas omitting the kinetic HoTs reduces the macroscopic subproblem to the reference mixture NS-DG solver.

\section{Numerical tests}\label{sec:4}

Three benchmark problems are considered to assess the accuracy, spatial resolution, asymptotic preserving property, and computational efficiency of GSIS-DG. The normal shock wave examines species-dependent nonequilibrium effects and high-order spatial accuracy, the lid-driven cavity covers flow regimes from highly rarefied to continuum, and the nozzle flow evaluates the performance of the method for an expanding mixture.
The species mole fraction is $\chi_s=n_s/(n_1+n_2)$, and the number-density ratio is $\beta_n=n_2/n_1$. The normal-shock and nozzle calculations use a Maxwell mixture with $\beta_m=m_2/m_1=10$ and $\beta_d=d_2/d_1=1$. 
For the lid-driven cavity, $\beta_m=\beta_n=1$, so that the mixture model reduces to the monatomic single-gas limit.

The convergence of both CIS-DG and GSIS-DG is monitored using the volume-weighted relative change of the macroscopic moments,
\begin{equation}
E^k=\max_{s,\,Q\in\{\rho_s,\bm u_s,T_s\}}
\frac{\left[\sum_j\omega_j\|Q_j^k-Q_j^{k-1}\|_2^2\right]^{1/2}}
{\left[\sum_j\omega_j\|Q_j^{k-1}\|_2^2\right]^{1/2}},
\label{eq:dg-convergence}
\end{equation}
where $j$ indexes physical cubature points and $\omega_j$ includes the quadrature and element-volume weights. The norm is the absolute value for scalars and the Euclidean norm for velocity, and the maximum is taken over the three monitored quantities and both species. 
The prescribed tolerance is $\epsilon_{\mathrm{tol}}=10^{-5}$ for the normal shock and $\epsilon_{\mathrm{tol}}=10^{-6}$ for the cavity and nozzle flows. In the following, iteration counts are reported as the number of kinetic updates; each GSIS-DG outer iteration contains one kinetic update followed by the inner iterations of the macroscopic synthetic equations.

To quantify the spatial accuracy, a normalized weighted $L_2$ error is used throughout the following tests,
\begin{equation}
E_2(Q)=
\frac{
\left[
\sum_j \omega_j
\left\|Q_j-Q_j^{\mathrm{ref}}\right\|_2^2
\right]^{1/2}
}{
\mathcal{N}_Q
},
\label{eq:dg-l2-error}
\end{equation}
where $Q^{\mathrm{ref}}$ denotes the corresponding reference solution, $\omega_j$ is the spatial weight, and $\mathcal{N}_Q$ is the normalization specified for each comparison. The same definition is applied to either scalar or vector macroscopic quantities. For relative errors, $\mathcal N_Q=[\sum_j\omega_j\|Q_j^{\mathrm{ref}}\|_2^2]^{1/2}$. For errors normalized by a fixed physical scale $S_Q$, $\mathcal N_Q=S_Q(\sum_j\omega_j)^{1/2}$. 

\subsection{Normal shock wave}\label{sec:dg-shock}

\begin{figure}[!ht]
\centering
\includegraphics[scale=0.85,clip = true]{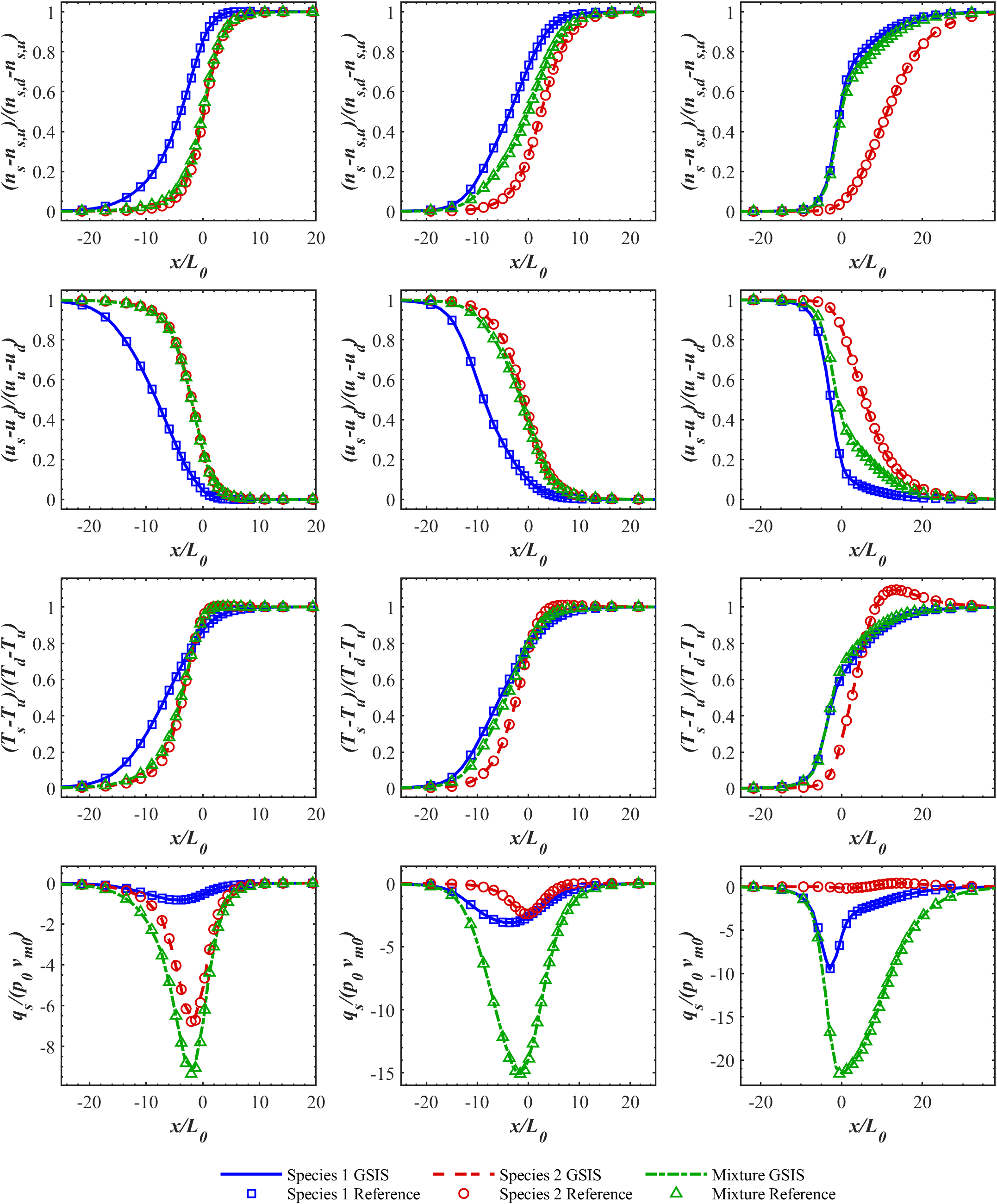}
\caption{Normal-shock profiles obtained by GSIS-DG (lines) and DSMC (symbols) at $\mathrm{Ma}=5$. Columns correspond to $\chi_1=0.1$, 0.5, and 0.9; rows show the normalized number density, streamwise velocity, temperature, and heat flux. Blue squares, red circles, and green triangles identify species 1, species 2, and the mixture, respectively. The heat-flux scale is $p_0v_{m0}$, where $p_0=n_uk_BT_u$ and $v_{m0}=\sqrt{2k_BT_u/m_1}=\sqrt{2}\,v_0$. The calculation and DSMC profiles are aligned independently using the midpoint of the normalized total number density; the horizontal coordinate $x/L_0$ is reported relative to this aligned shock center.}
\label{fig:dg-shock-profiles}
\end{figure}

A stationary normal shock wave is first considered to assess the accuracy, spatial convergence, and iterative convergence of GSIS-DG without complications from gas-surface boundary conditions. The molecular mass of the light species, the density of the number of the upstream mixture, and the temperature upstream are taken as reference quantities, i.e. $m_0=m_1$, $n_0=n_u$, and $T_0=T_u$. The characteristic length $L_0$ is the upstream mean free path of the light species, such that $\Kn_1=1$. A Maxwell mixture with $\beta_m=10$ and $\beta_d=d_2/d_1=1$ is considered at $\mathrm{Ma}=5$, where the Mach number is defined using the upstream mixture sound speed $v_s=\sqrt{5k_BT_u/(3m_{\rm mix})}$, with $m_{\rm mix}=\sum_s m_s\chi_s$. Three upstream compositions, $\chi_1=0.1$, $0.5$, and $0.9$, are examined. The computational domain $[0,100L_0]$ is discretized by 24 uniform elements with $p=2$. Separate molecular-velocity grids are employed for the two species to accommodate their different thermal speeds, with the velocity ranges and resolutions adjusted slightly with the mixture composition. The upstream and downstream equilibrium states are related by the classical Rankine--Hugoniot conditions and are used to initialize the two halves of the computational domain.


Figure~\ref{fig:dg-shock-profiles} compares the GSIS-DG solutions with DSMC for the three mixture compositions. The cell width is $\Delta x/L_0\simeq4.2$, so that a single cell spans more than four upstream mean free paths of the light species. With $p=2$, only 72 spatial degrees of freedom are used per species and molecular velocity. Nevertheless, the number density, streamwise velocity, temperature, and heat flux agree closely with DSMC over all three compositions. The characteristic species-dependent nonequilibrium behavior is also reproduced: the density and velocity transitions of the light species occur upstream of those of the heavy species, and a pronounced heavy-species temperature overshoot appears at $\chi_1=0.9$. The largest visible discrepancy is found in the light-species heat flux at $\chi_1=0.9$, where the difference is approximately $3.5\%$ of the DSMC peak value.

\begin{table}[!t]
\centering
\caption{
Relative $L_2$ errors of the normalized light-species velocity with respect to DSMC for the normal shock at $\mathrm{Ma}=5$, obtained with second-order GSIS-FV and GSIS--DG with $p=1,2,3$. The GSIS-FV scheme follows Ref.~\cite{zeng2024gsismixture}.
}
\label{tab:nsw-dsmc-accuracy}
\small
\setlength{\tabcolsep}{5pt}
\begin{tabular}{cr cr cr cr}
\toprule
\multicolumn{2}{c}{GSIS-FV} & \multicolumn{2}{c}{$p=1$} & \multicolumn{2}{c}{$p=2$} & \multicolumn{2}{c}{$p=3$}\\
\cmidrule(lr){1-2}\cmidrule(lr){3-4}\cmidrule(lr){5-6}\cmidrule(lr){7-8}
Mesh & Error (\%) & Mesh & Error (\%) & Mesh & Error (\%) & Mesh & Error (\%)\\ \midrule
18 & 10.26 & 18 & 2.79 & 12 & 0.81 & 12 & 0.27\\
24 & 6.28 & 24 & 1.67 & 18 & 0.31 & 18 & 0.13\\
32 & 3.91 & 32 & 0.95 & 24 & 0.18 & 24 & 0.12\\
48 & 1.84 & 48 & 0.44 & 32 & 0.13 & 32 & 0.12\\
64 & 0.73 & 64 & 0.27 & 48 & 0.12 & 48 & 0.12\\
\bottomrule
\end{tabular}
\end{table}

\begin{figure}[!htbp]
\centering
\includegraphics[width=0.98\textwidth]{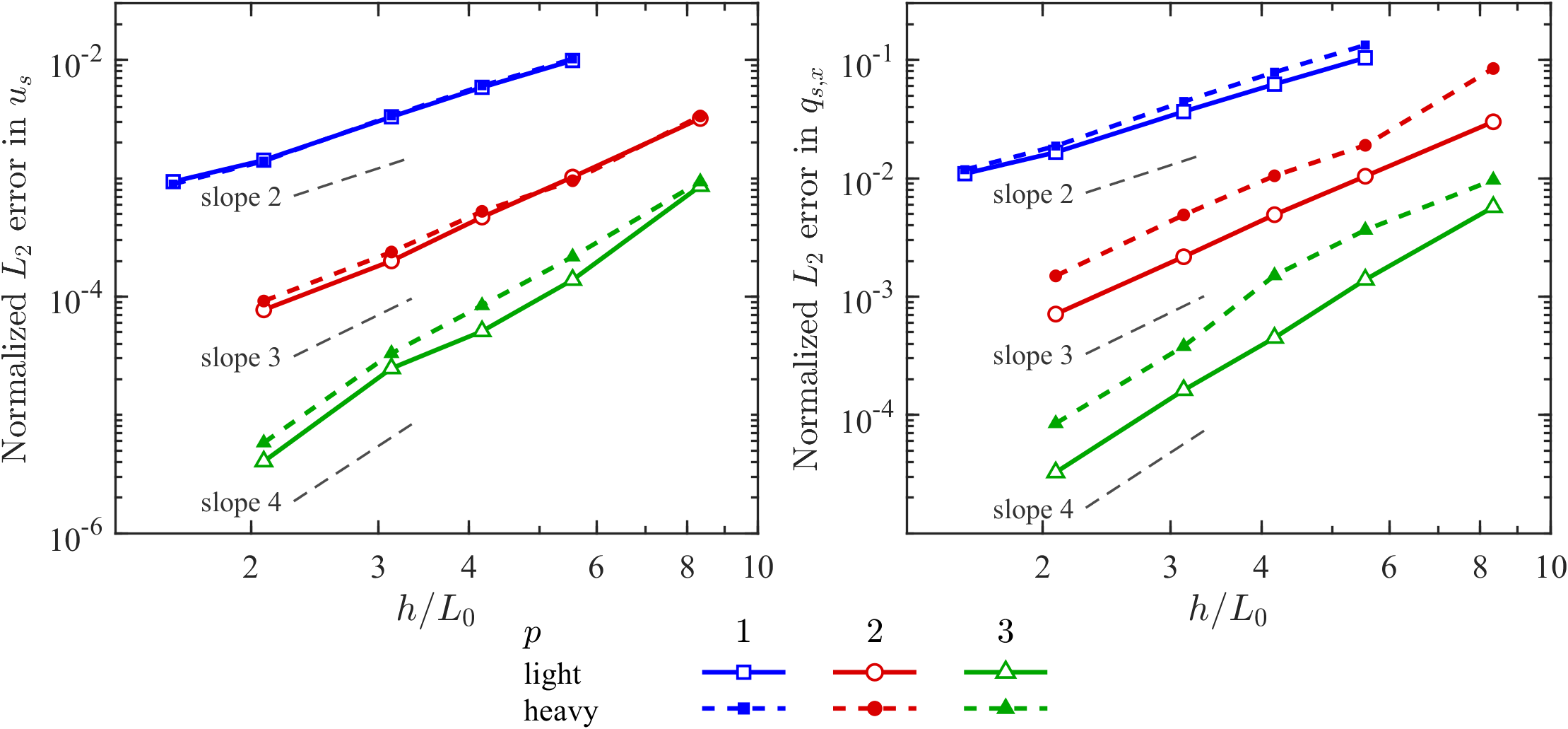}
\caption{
Spatial convergence of the species velocity (left) and heat flux (right) for the equimolar normal shock at $\mathrm{Ma}=5$, using $p=1,2,3$. Errors are evaluated against same-degree 96-element reference solutions. Solid lines with open symbols and dashed lines with filled symbols denote the light and heavy species, respectively; gray dashed guides indicate the nominal slopes $2$--$4$.
}
\label{fig:nsw-spatial-order}
\end{figure}

Table~\ref{tab:nsw-dsmc-accuracy} compares the relative errors obtained with different polynomial degrees, using the DSMC solution as the reference. Since the $p=1$ DG discretization is formally second-order accurate, results from the second-order finite-volume GSIS are also included for a same-order comparison. On the same physical mesh, GSIS-DG consistently yields smaller errors than GSIS-FV; moreover, the $p=1$ DG solution with 24 elements is already more accurate than the finite-volume solution with 48 cells. Increasing the polynomial degree further reduces the spatial resolution required for a given accuracy. In particular, the $p=3$ solution with only 12 elements reaches approximately the same accuracy as the $p=1$ solution with 64 elements. To further examine the spatial convergence order of the DG discretization, 96-element solutions with the corresponding polynomial degree are used as numerical references. The velocity and heat-flux errors are evaluated using Eq.~\eqref{eq:dg-l2-error}, with $S_u=u_u-u_d$ and $S_q=p_0v_{m0}$ as the normalization scales for velocity and heat flux, respectively. Figure~\ref{fig:nsw-spatial-order} shows that the observed convergence rates for $p=1,2,3$ are close to the nominal $p+1$ order for both species, while the heat-flux errors exhibit the same systematic improvement with increasing polynomial degree.

\begin{figure}[!htbp]
\centering
\includegraphics[width=0.6\textwidth,keepaspectratio]{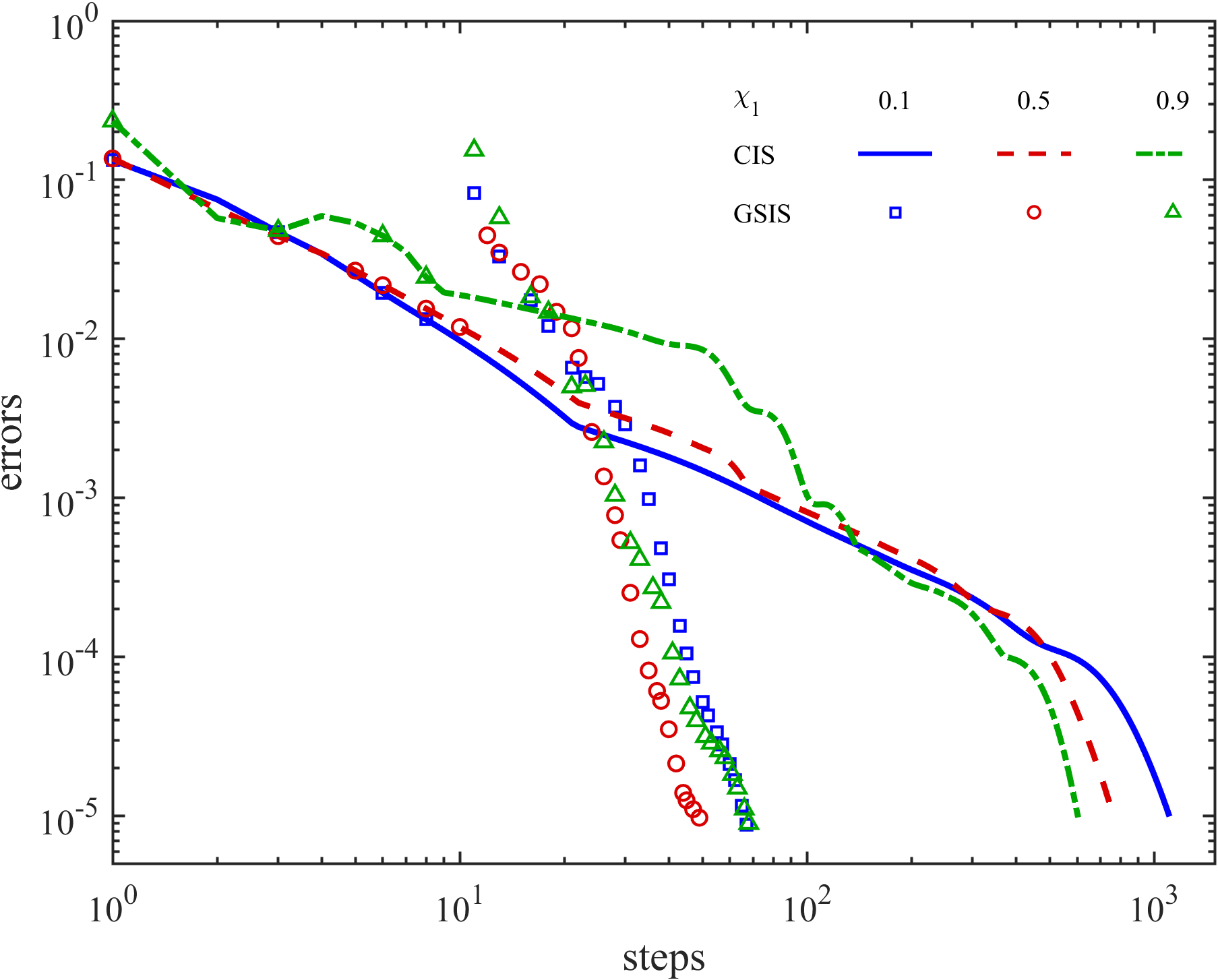}
\caption{
Convergence histories for the normal shock at $\mathrm{Ma}=5$. Blue, red, and green denote $\chi_1=0.1$, 0.5, and 0.9, respectively. Solid, dashed, and dash-dotted lines show CIS-DG; matching open symbols show GSIS-DG. All calculations use 24 physical elements with $p=2$, matched initial states, and identical molecular-velocity discretizations within each composition. The horizontal axis counts kinetic updates; for GSIS-DG, each outer iteration contains one kinetic update followed by the macroscopic synthetic iterations. GSIS-DG is initialized by ten CIS kinetic updates before the synthetic iteration is activated.
}
\label{fig:dg-shock-convergence}
\end{figure}

Figure~\ref{fig:dg-shock-convergence} compares the convergence histories of CIS-DG and GSIS-DG using identical initial states and discretizations. GSIS-DG converges within several tens of kinetic updates for all three compositions, whereas CIS-DG requires several hundred to more than one thousand updates.
The convergence rate of CIS-DG varies noticeably with the mixture composition, while that of GSIS-DG remains within a relatively narrow range because the macroscopic synthetic equations transfer information across the entire computational domain rather than relying solely on molecular streaming.
All simulations are performed on the same AMD EPYC 7742 processor using eight CPU cores. For the equimolar case, the wall-clock time is reduced from 259.7 to 27.7~s, corresponding to approximately one order of magnitude acceleration.

\subsection{Lid-driven cavity flow}
\label{sec:dg-cavity}

The lid-driven cavity flow is considered to assess GSIS-DG from the highly rarefied to the continuum regime. The equal-species configuration has $\beta_m=\beta_n=1$, for which the mixture model recovers the monatomic limit~\cite{li2024kinetic}. The square cavity has side length $L_0$ and is initially filled with an equilibrium gas at $\rho=1$, $\bm{u}=0$, and $T=T_0=1$. The upper wall moves with $u_w=0.2066v_0$, corresponding to $\mathrm{Ma}=0.16$, while the other walls remain stationary. All walls are maintained at $T_w=T_0$ and employ fully diffuse reflection. The cases at $\Kn_1=10$ and 0.075 use a $4\times4$ physical mesh with $p=3$ and one wall-adjacent layer of height $0.05L_0$, whereas the $\mathrm{Re}=1000$ case uses a $6\times6$ mesh with $p=4$ and a first-layer height of $0.02L_0$. A uniform $60\times60$ molecular-velocity grid over $[-6,6]^2$ is used in all cases.

\begin{figure}[!htbp]
\centering
\begin{minipage}[c]{0.44\textwidth}\centering
\includegraphics[width=\linewidth]{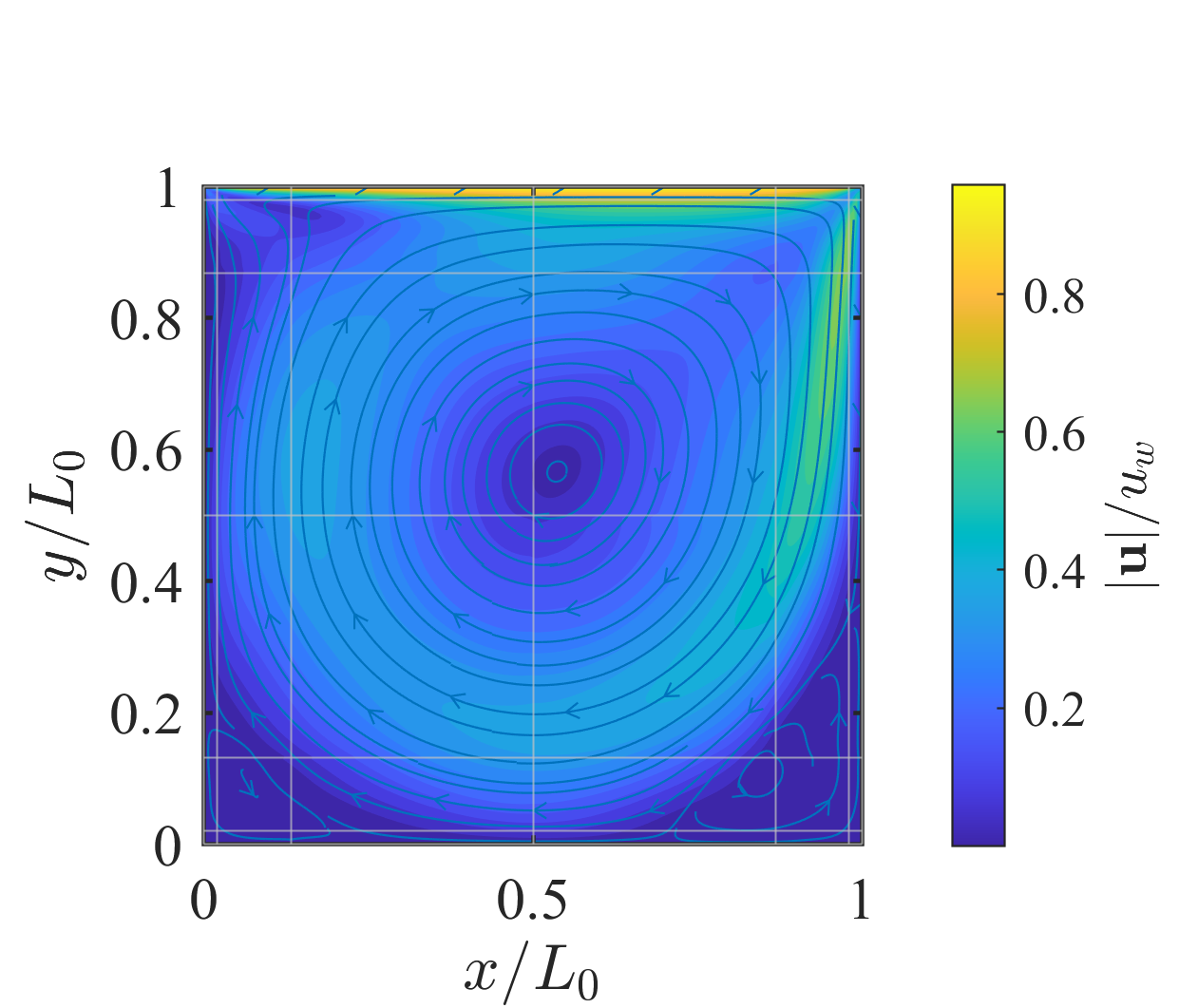}
\end{minipage}\hspace{0.01\textwidth}%
\begin{minipage}[c]{0.44\textwidth}\centering
\includegraphics[width=\linewidth]{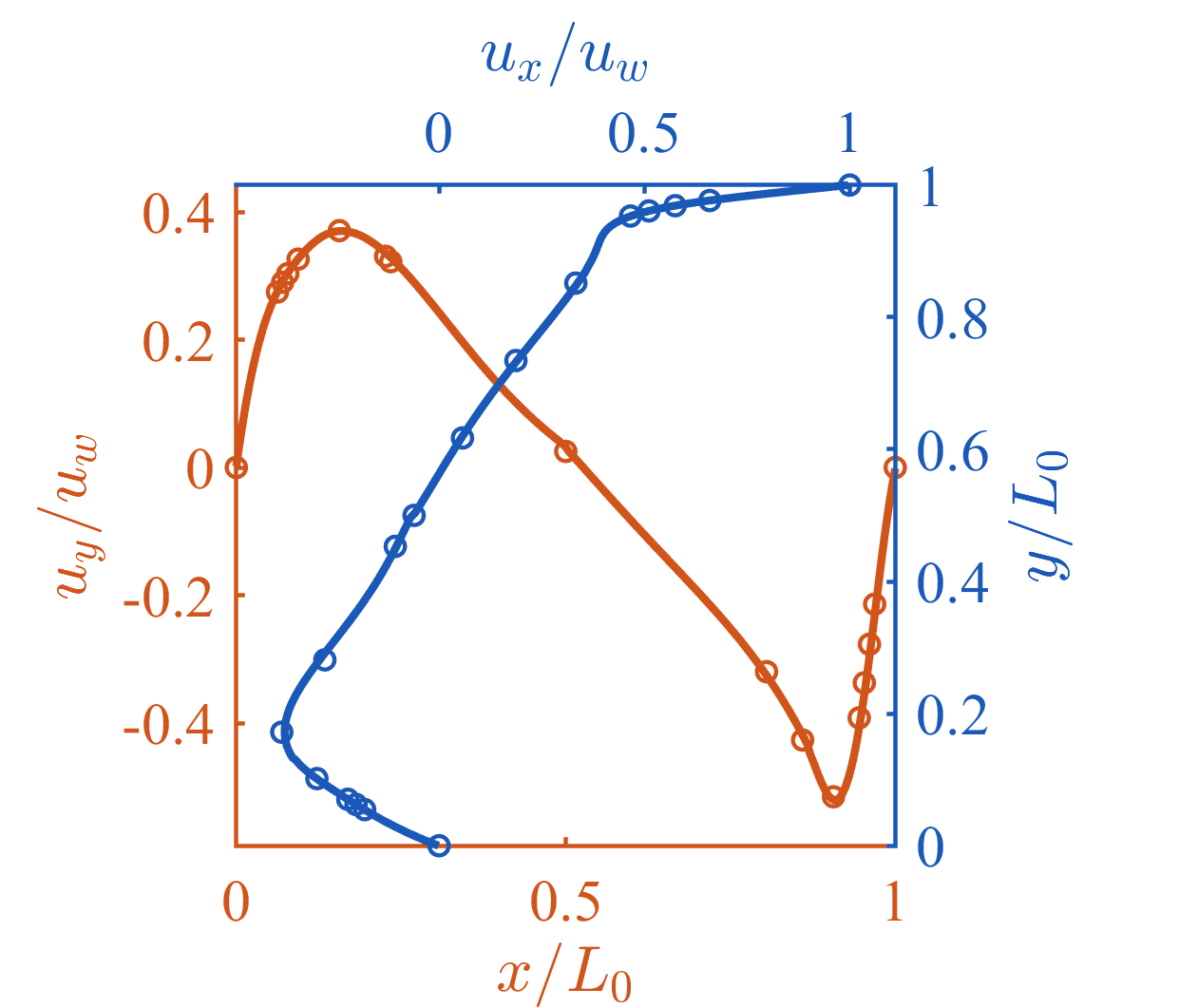}
\end{minipage}\par\vspace{0.05em}
{\small\bfseries (a) $\mathrm{Re}=1000$}\par\vspace{0.35em}
\begin{minipage}[c]{0.44\textwidth}\centering
\includegraphics[width=\linewidth]{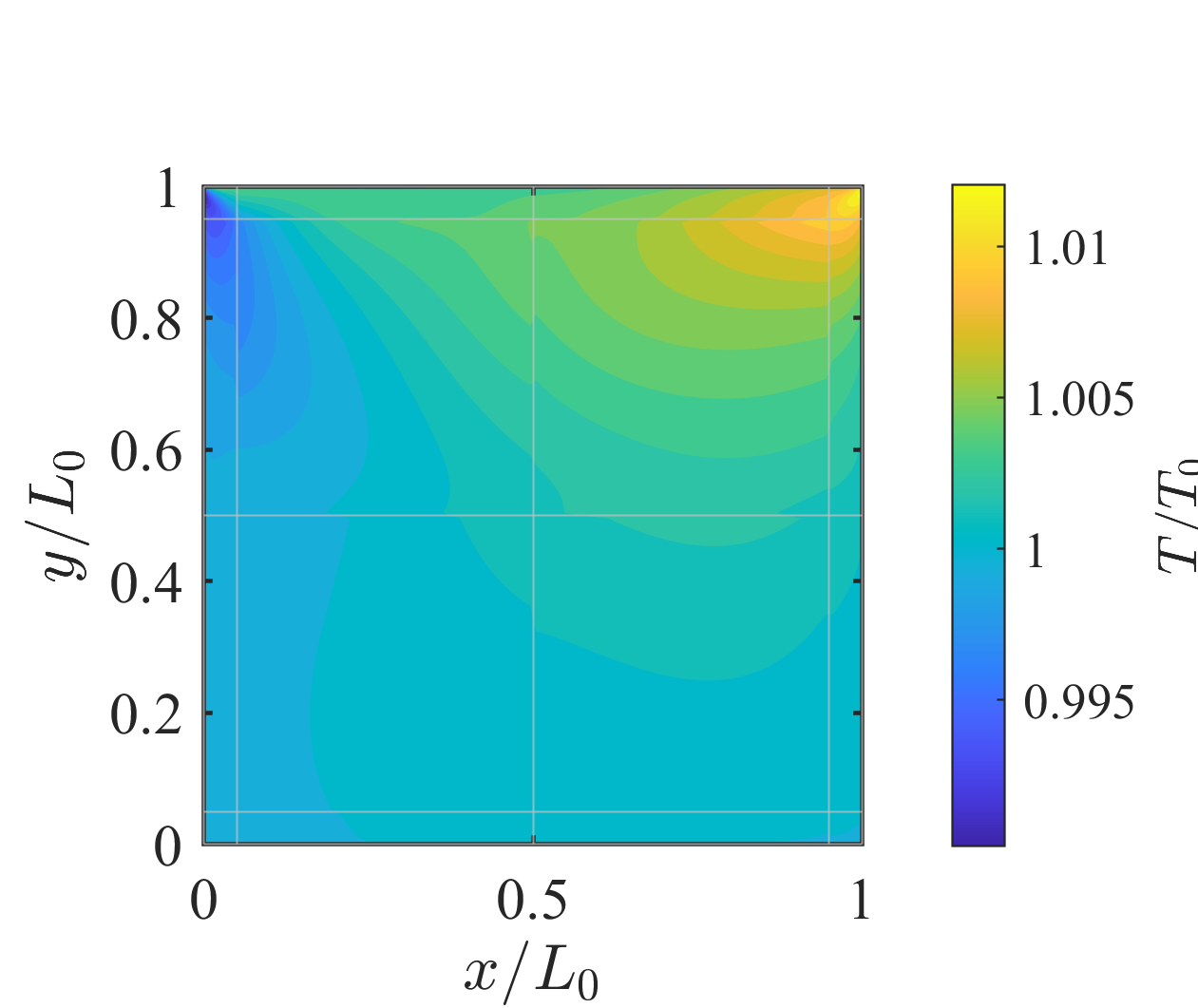}
\end{minipage}\hspace{0.01\textwidth}%
\begin{minipage}[c]{0.44\textwidth}\centering
\includegraphics[width=\linewidth]{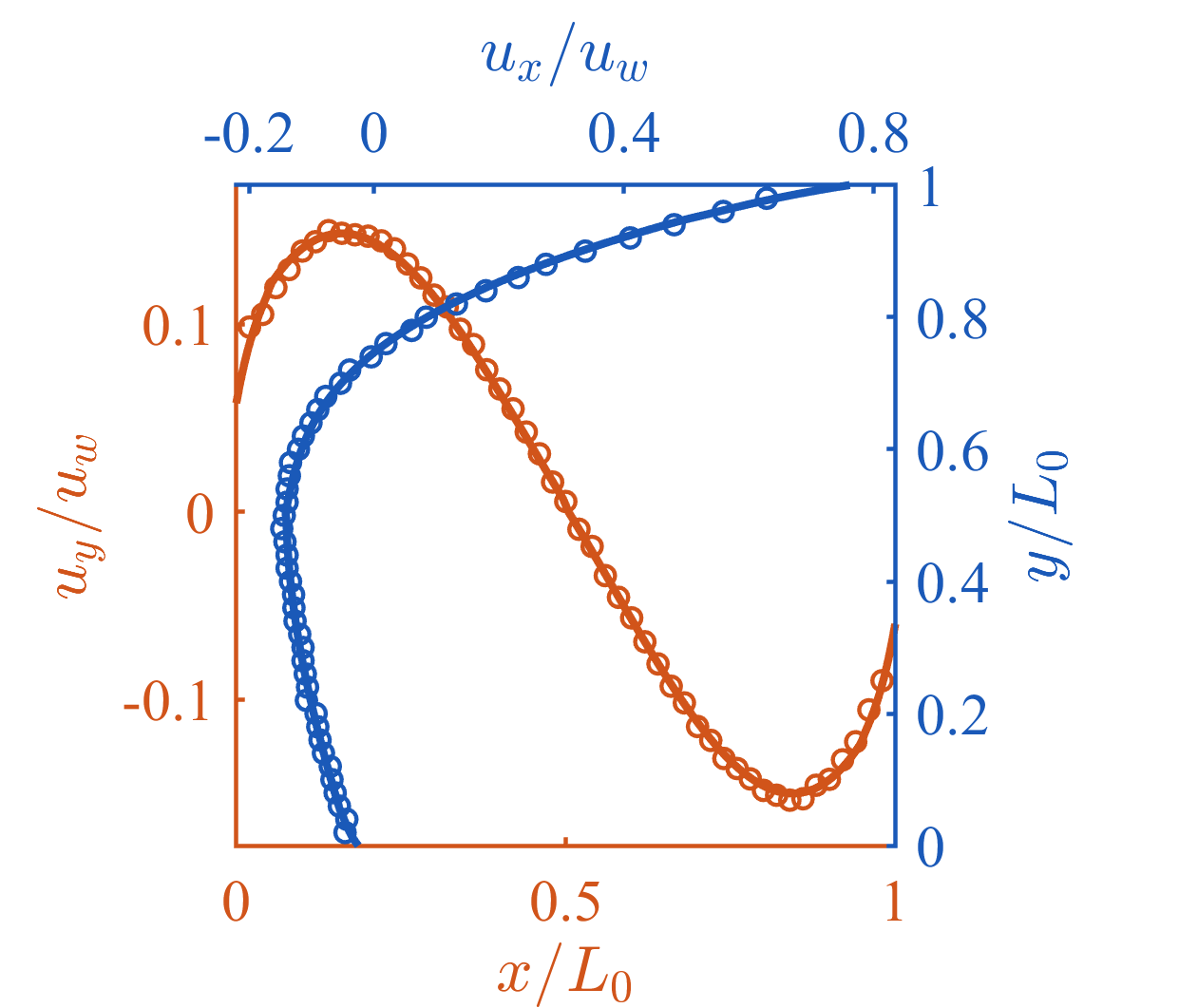}
\end{minipage}\par\vspace{0.05em}
{\small\bfseries (b) $\Kn_1=0.075$}\par\vspace{0.35em}
\begin{minipage}[c]{0.44\textwidth}\centering
\includegraphics[width=\linewidth]{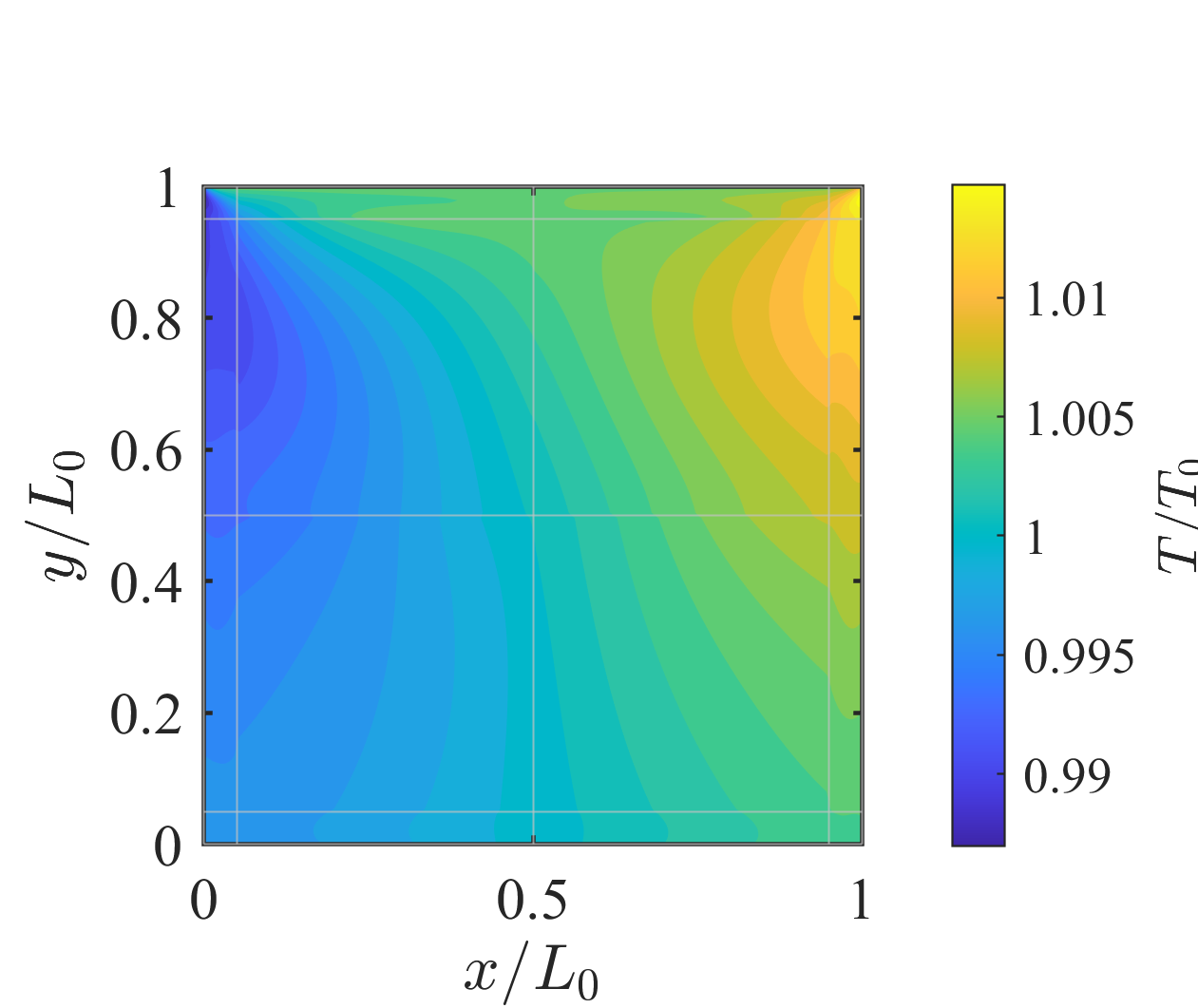}
\end{minipage}\hspace{0.01\textwidth}%
\begin{minipage}[c]{0.44\textwidth}\centering
\includegraphics[width=\linewidth]{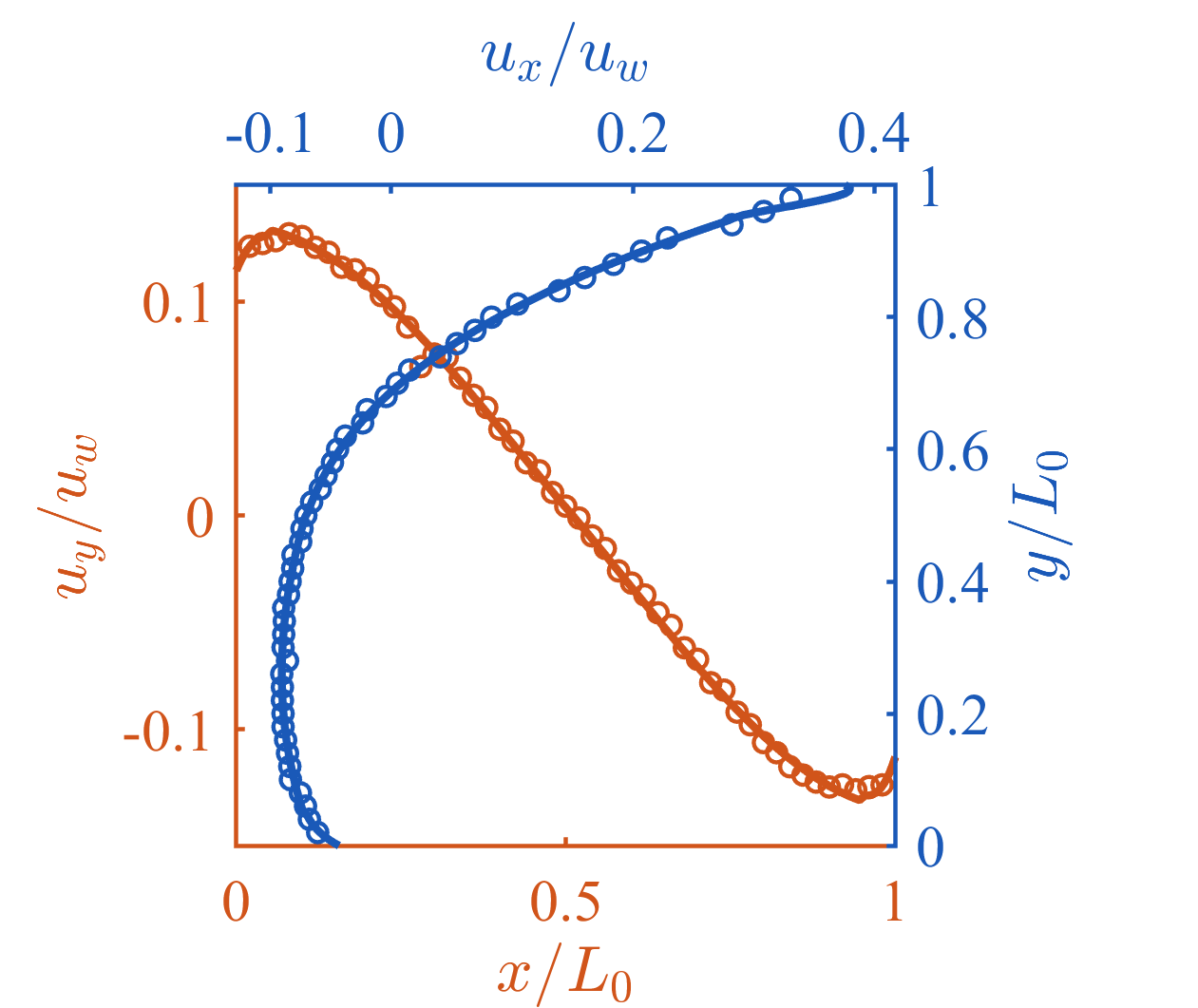}
\end{minipage}\par\vspace{0.05em}
{\small\bfseries (c) $\Kn_1=10$}
\caption{
Lid-driven cavity at (a) $\mathrm{Re}=1000$, (b) $\Kn_1=0.075$, and (c) $\Kn_1=10$. Left: GSIS-DG flow fields. Right: centerline velocities; solid lines are GSIS-DG and circles are Ghia \emph{et al.}~\cite{ghia1982highre} in (a) or DUGKS~\cite{yuan2021novel} in (b,c).
}
\label{fig:dg-cavity-solutions}
\end{figure}

Figure~\ref{fig:dg-cavity-solutions} compares the GSIS-DG solutions with the benchmark data of Ghia \emph{et al.}~\cite{ghia1982highre} at $\mathrm{Re}=1000$ and the DUGKS results of Yuan \emph{et al.}~\cite{yuan2021novel} at $\Kn_1=10$ and 0.075. Good agreement is obtained for the centerline velocity profiles over all three flow regimes. At $\mathrm{Re}=1000$, the primary vortex and the two secondary vortices near the lower corners are reproduced using only $6\times6$ physical elements. Since the element size is then orders of magnitude larger than the molecular mean free path, the recovery of the continuum benchmark confirms that the asymptotic-preserving property is retained under the high-order DG discretization. The two rarefied cases are also accurately captured on the $4\times4$ mesh, including the characteristic heating and cooling near the moving-lid corners.


Table~\ref{tab:dg-cavity-cost} compares the convergence performance of CIS-DG and GSIS-DG. At $\Kn_1=10$, both methods converge in 33 kinetic updates, with wall-clock times of 10.3 and 10.7~s, respectively. The synthetic equations therefore provide no reduction in iteration count and incur a small additional cost in the highly rarefied regime. At $\Kn_1=0.075$, GSIS-DG reduces the kinetic updates from 281 to 48 and the wall-clock time from 82.1 to 15.7~s, corresponding to speedup factors of 5.9 and 5.2, respectively. The advantage becomes much more pronounced toward the continuum regime: at $\mathrm{Re}=1000$, GSIS-DG converges in 62 kinetic updates and 256.7~s, whereas CIS-DG remains unconverged after $10^4$ updates, giving an iteration-count speedup greater than 161. This trend reflects the increasing importance of the macroscopic synthetic equations in transferring information over the entire physical domain as the Knudsen number decreases.

\begin{table}[!htbp]
\centering
\caption{
Computational performance of CIS-DG and GSIS-DG for the lid-driven cavity. The cases at $\Kn_1=10$ and 0.075 use a $4\times4$ mesh with $p=3$, and the $\mathrm{Re}=1000$ case uses a $6\times6$ mesh with $p=4$.
}
\label{tab:dg-cavity-cost}
\begin{tabular}{c cc cc cc}
\toprule
& \multicolumn{2}{c}{CIS-DG} & \multicolumn{2}{c}{GSIS-DG} & \multicolumn{2}{c}{Speedup ratio} \\
\cmidrule(lr){2-3}\cmidrule(lr){4-5}\cmidrule(lr){6-7}
case & steps & wall time (s) & steps & wall time (s) & in steps & in time \\
\midrule
$\Kn_1=10$ & 33 & 10.3 & 33 & 10.7 & 1.0 & 0.97 \\
$\Kn_1=0.075$ & 281 & 82.1 & 48 & 15.7 & 5.9 & 5.2 \\
$\mathrm{Re}=1000$ & $>10{,}000$ & -- & 62 & 256.7 & $>161$ & -- \\
\bottomrule
\end{tabular}
\end{table}

To directly compare the spatial resolution of DG and finite-volume discretizations at the same formal order of accuracy, an additional study is performed at $\Kn_1=0.075$ using second-order GSIS-DG and GSIS-FV on uniform square meshes. The GSIS-DG discretization employs $Q_1$ polynomials with four spatial degrees of freedom per element, whereas GSIS-FV stores one cell-average value per cell; both use the same $60\times60$ molecular-velocity grid. The $128\times128$ GSIS-FV solution is taken as the reference, and the relative cell-average velocity error is evaluated using Eq.~\eqref{eq:dg-l2-error}, with $Q_j=\overline{\bm{u}}_j$, cell-area weights $\omega_j$, and the reference-solution norm for $\mathcal N_Q$. As shown in Table~\ref{tab:dg-fv-cavity-accuracy}, for the same spatial degrees of freedom, the $16\times16$ GSIS-DG solution gives a relative error of $0.7\%$, compared with $2.6\%$ for the $32\times32$ GSIS-FV solution. Moreover, the $16\times16$ DG solution is more accurate than the $64\times64$ FV solution, whose error is $1.1\%$, while using only one quarter of the spatial degrees of freedom. The wall-clock time decreases from $96.7$s for FV to $77.3$s for DG. On the same physical mesh, however, GSIS-DG generally requires more iterations than GSIS-FV because of the smaller admissible macroscopic CFL number in the present DG implementation, which motivates future hybrid FV-DG and $hp$-adaptive strategies.


\begin{table}[!htbp]
\centering
\caption{Spatial accuracy and recorded computational cost of second-order GSIS-DG and GSIS-FV for the cavity at $\Kn_1=0.075$. The reference is the $128\times128$ GSIS-FV solution.}
\label{tab:dg-fv-cavity-accuracy}
\begingroup
\setlength{\tabcolsep}{5pt}
\begin{tabular}{crrrrrr}
\toprule
& \multicolumn{3}{c}{GSIS-DG} & \multicolumn{3}{c}{GSIS-FV} \\
\cmidrule(lr){2-4}\cmidrule(lr){5-7}
Mesh & Steps & Time (s) & $E_2$ (\%) & Steps & Time (s) & $E_2$ (\%) \\
\midrule
$4\times4$ & 37 & 1.2 & 5.7 & 26 & 0.1 & 22.5 \\
$8\times8$ & 65 & 10.8 & 1.8 & 28 & 0.3 & 10.0 \\
$16\times16$ & 100 & 77.3 & 0.7 & 35 & 1.6 & 4.8 \\
$32\times32$ & 227 & 888.2 & 0.5 & 60 & 12.8 & 2.6 \\
$64\times64$ & -- & -- & -- & 109 & 96.7 & 1.1 \\
\bottomrule
\end{tabular}
\endgroup
\end{table}

\subsection{Three-dimensional nozzle flow}\label{sec:dg-nozzle}

A three-dimensional converging--diverging nozzle is considered to assess the accuracy and computational efficiency of GSIS-DG for expanding binary-mixture flows. The nozzle has a length of $3L_0$, where $L_0$ is the inlet diameter. Its meridional profile passes through $A(0,0.5)$, $B(0.5,0.5)$, $C(1,0.25)$, and $D(3,1)$, giving throat and outlet diameters of $0.5L_0$ and $2L_0$, respectively. Although the geometry is axisymmetric, all GSIS-DG calculations are performed on the full three-dimensional domain. The equimolar Maxwell mixture considered above is used, with $\beta_d=d_2/d_1=1$. At the inlet, both species follow Maxwellian velocity distributions at $T_0$ with $\mathrm{Ma}_{\mathrm{in}}=0.05$. The nozzle wall is maintained at $T_w=T_0/2$ with fully diffuse reflection, while the gas expands into vacuum at the outlet. Two rarefaction levels, $\Kn_1=0.1$ and 0.01, are considered.

\begin{figure}[!htbp]
\centering
\subfloat[Physical mesh]{\includegraphics[width=0.35\textwidth]{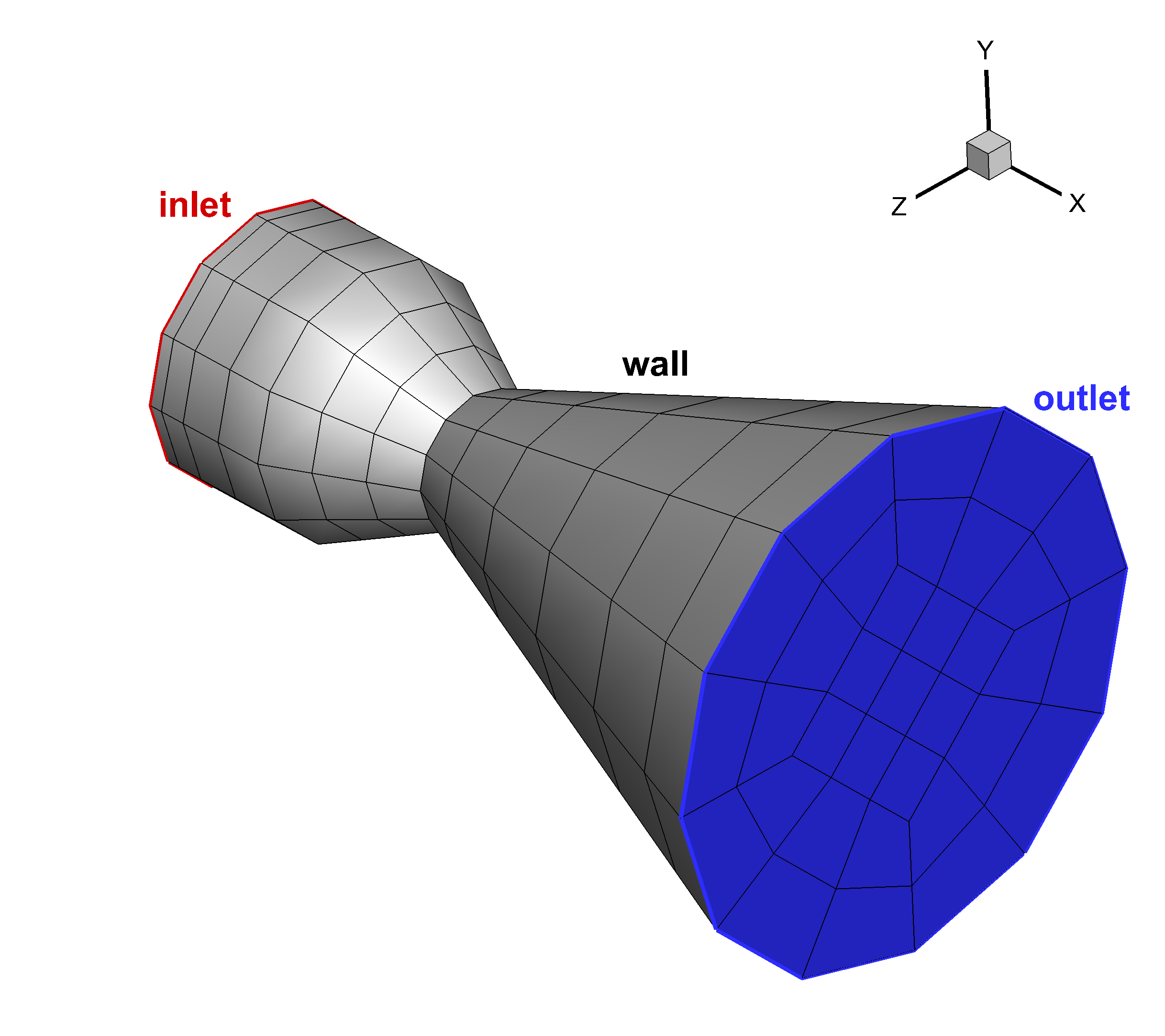}}\;
\subfloat[Meridional section]{\includegraphics[width=0.4\textwidth]{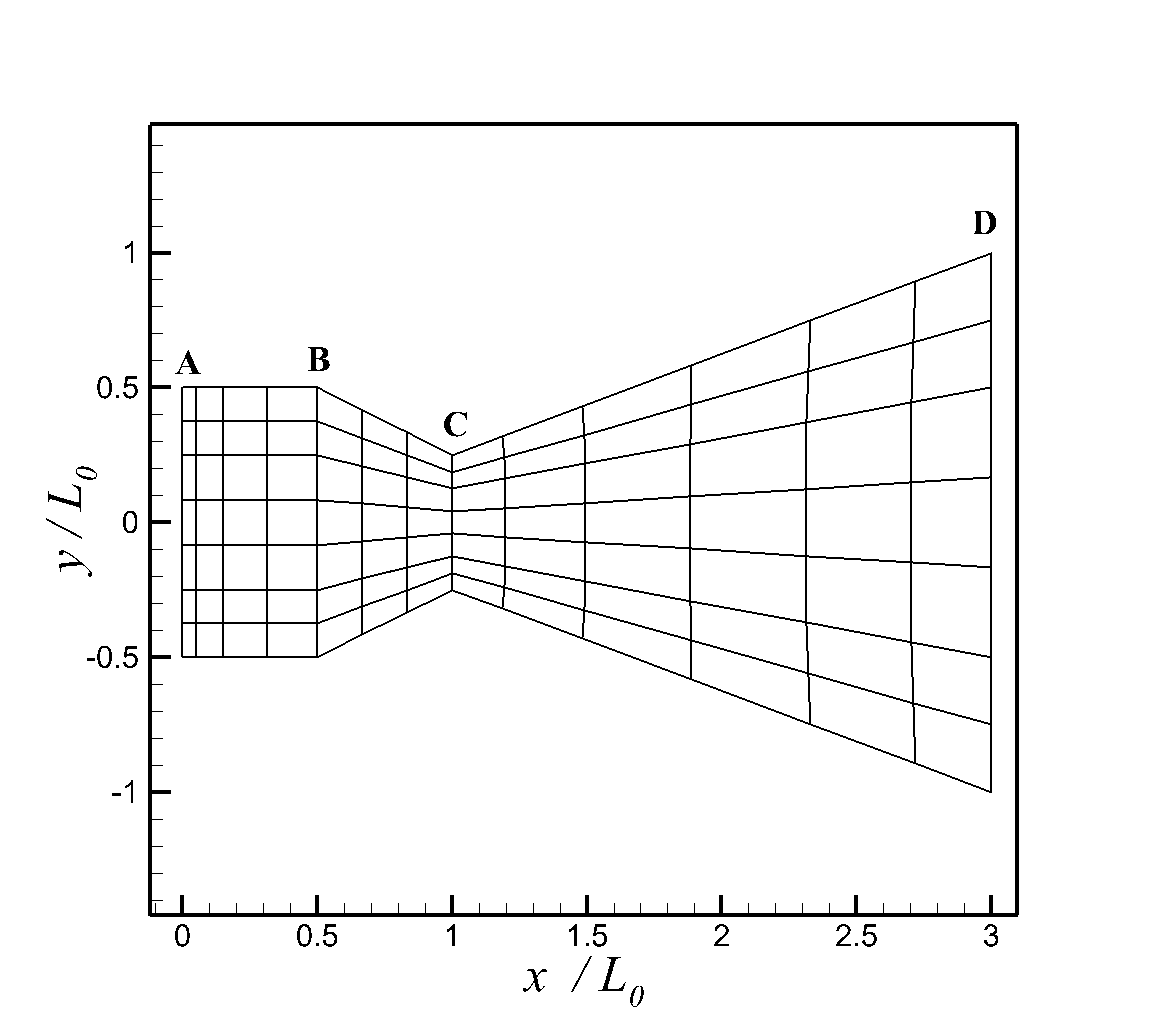}}\\[0.2em]
\subfloat[Light species velocity mesh]{\includegraphics[width=0.35\textwidth]{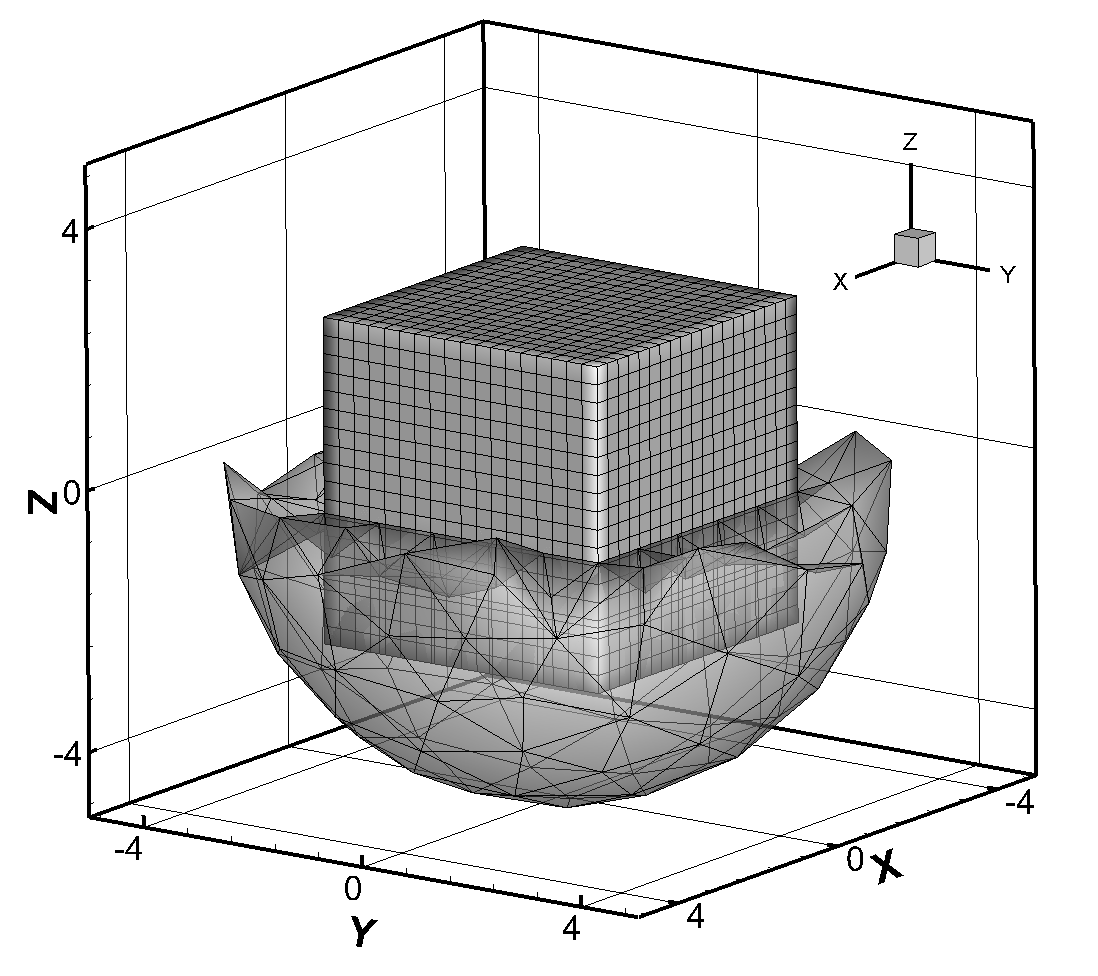}}\;
\subfloat[Heavy species velocity mesh]{\includegraphics[width=0.35\textwidth]{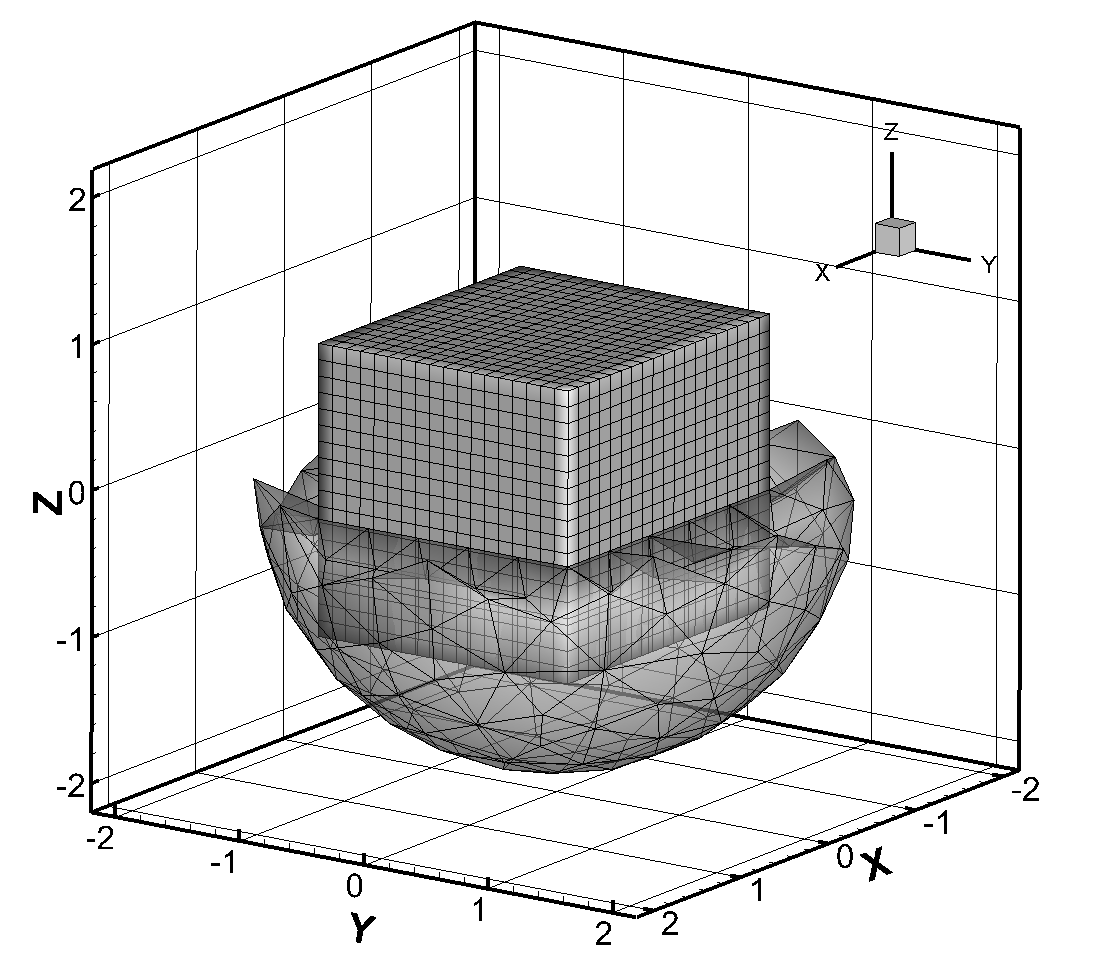}}
\caption{
Geometry and computational meshes for the nozzle flow. The upper row shows the full three-dimensional physical mesh and its meridional section; the physical domain contains 429 hexahedral elements with $p=2$. The lower row shows the hybrid molecular-velocity meshes for the light and heavy species, containing 8509 and 8625 velocity points, respectively. Molecular velocities are normalized by $v_0=\sqrt{k_BT_0/m_1}$.
}
\label{fig:dg-nozzle-mesh}
\end{figure}

The nozzle geometry and computational meshes are shown in Fig.~\ref{fig:dg-nozzle-mesh}. The physical domain is discretized by 429 hexahedral elements with $p=2$, corresponding to 27 nodal degrees of freedom in each element and 11,583 spatial degrees of freedom per species and molecular velocity. For comparison, the axisymmetric DSMC simulation at $\Kn_1=0.01$ in Ref.~\cite{zeng2024gsismixture} employed 373,761 physical cells. The hybrid structured--unstructured molecular-velocity discretization follows the previous GSIS study~\cite{zeng2024gsismixture}. The velocity grids contain 8509 and 8625 points for the light and heavy species, respectively. The same physical and molecular-velocity meshes are used at both Knudsen numbers.

\begin{figure}[!htbp]
\centering
\begin{minipage}[t]{0.485\textwidth}\centering
\includegraphics[width=\linewidth]{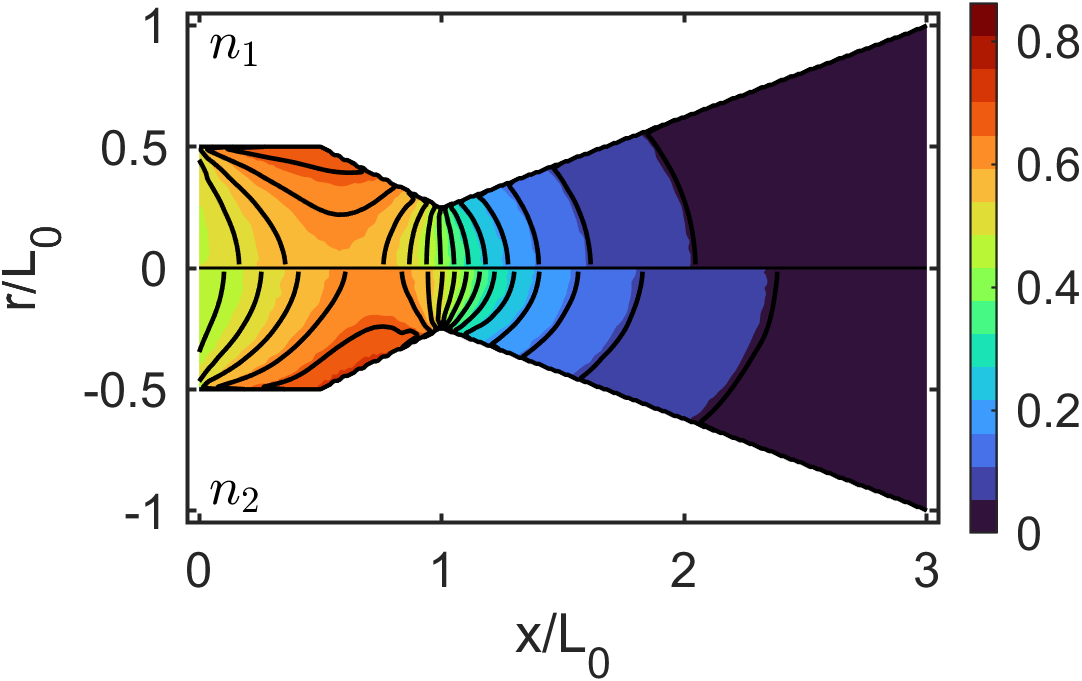}
\end{minipage}\hfill
\begin{minipage}[t]{0.485\textwidth}\centering
\includegraphics[width=\linewidth]{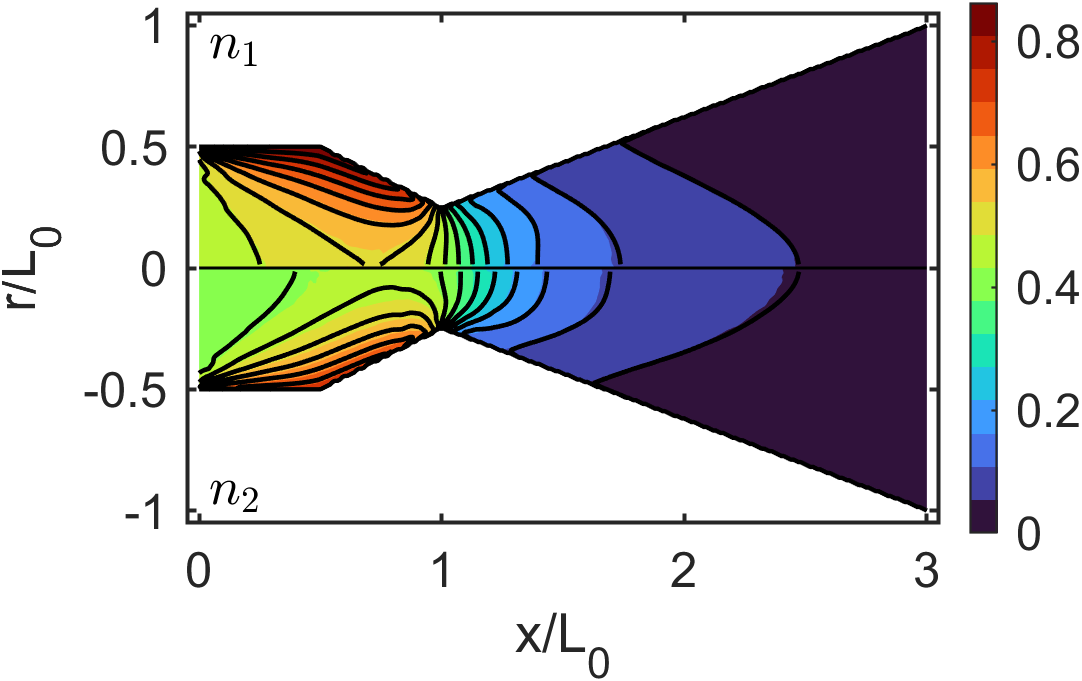}
\end{minipage}\\[0.35em]
\begin{minipage}[t]{0.485\textwidth}\centering
\includegraphics[width=\linewidth]{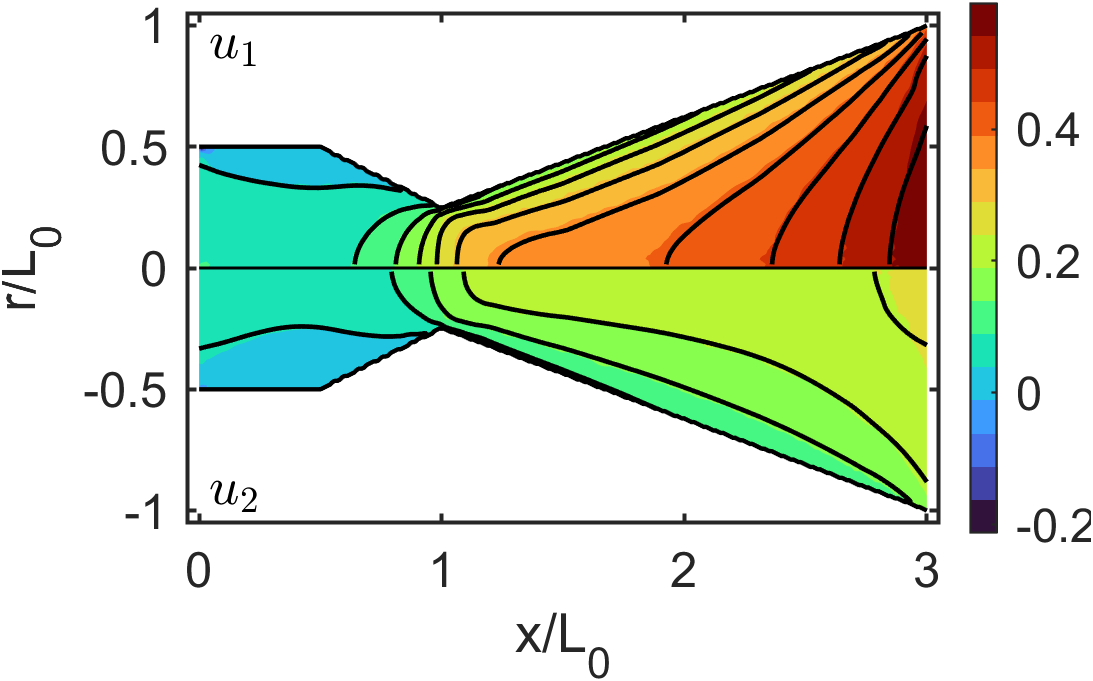}
\end{minipage}\hfill
\begin{minipage}[t]{0.485\textwidth}\centering
\includegraphics[width=\linewidth]{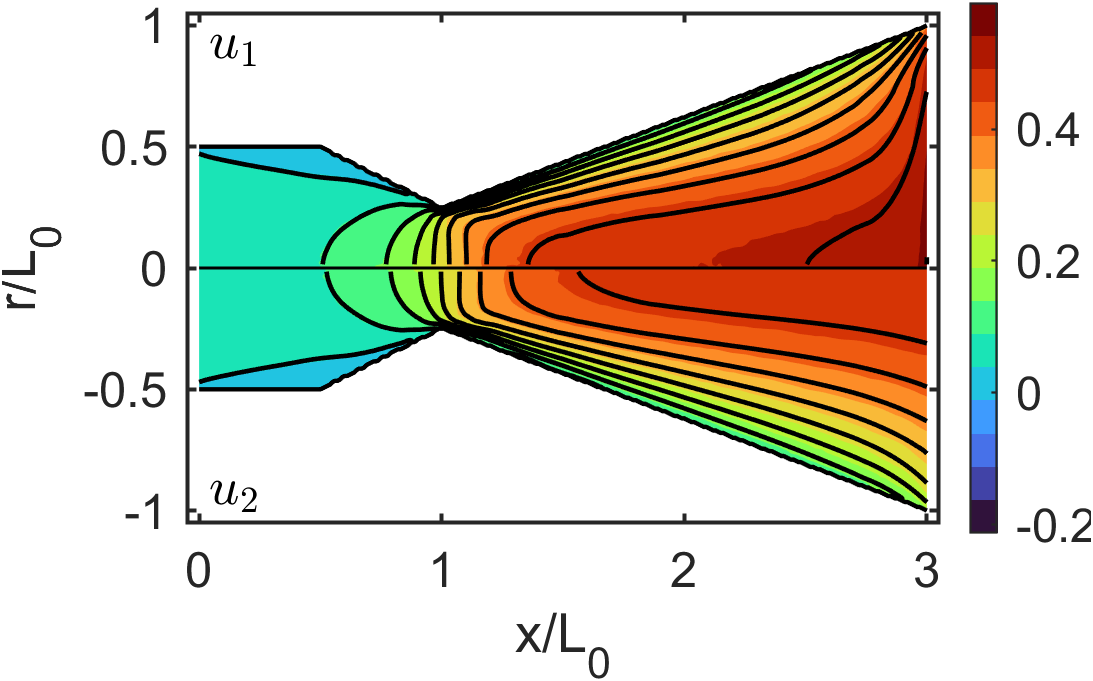}
\end{minipage}\\[0.35em]
\begin{minipage}[t]{0.485\textwidth}\centering
\includegraphics[width=\linewidth]{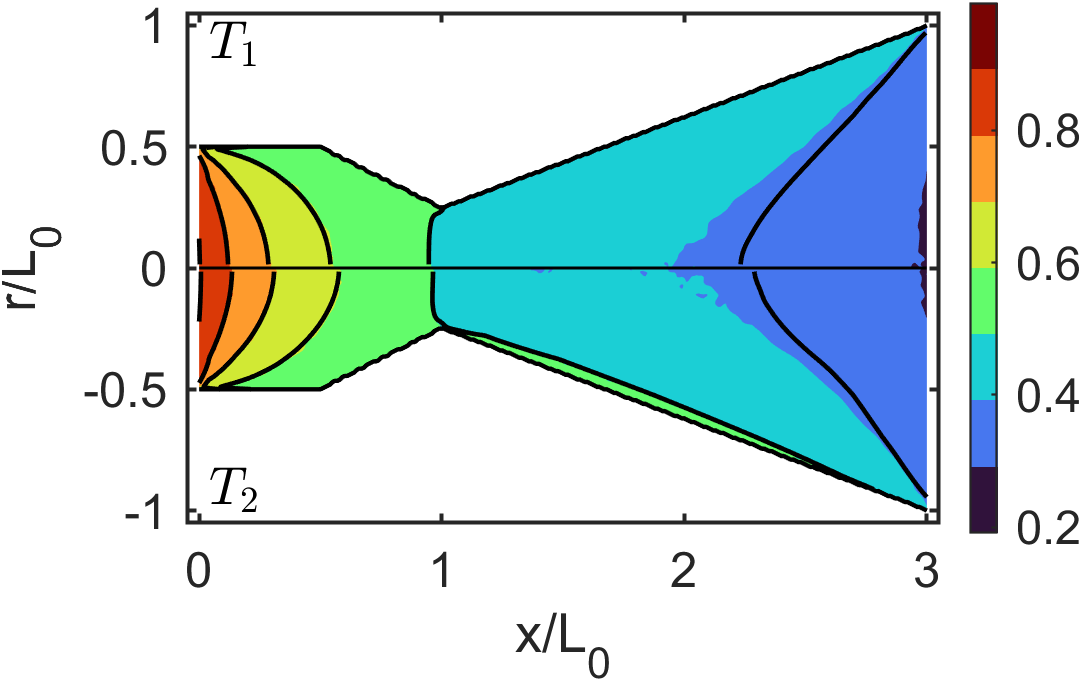}
\end{minipage}\hfill
\begin{minipage}[t]{0.485\textwidth}\centering
\includegraphics[width=\linewidth]{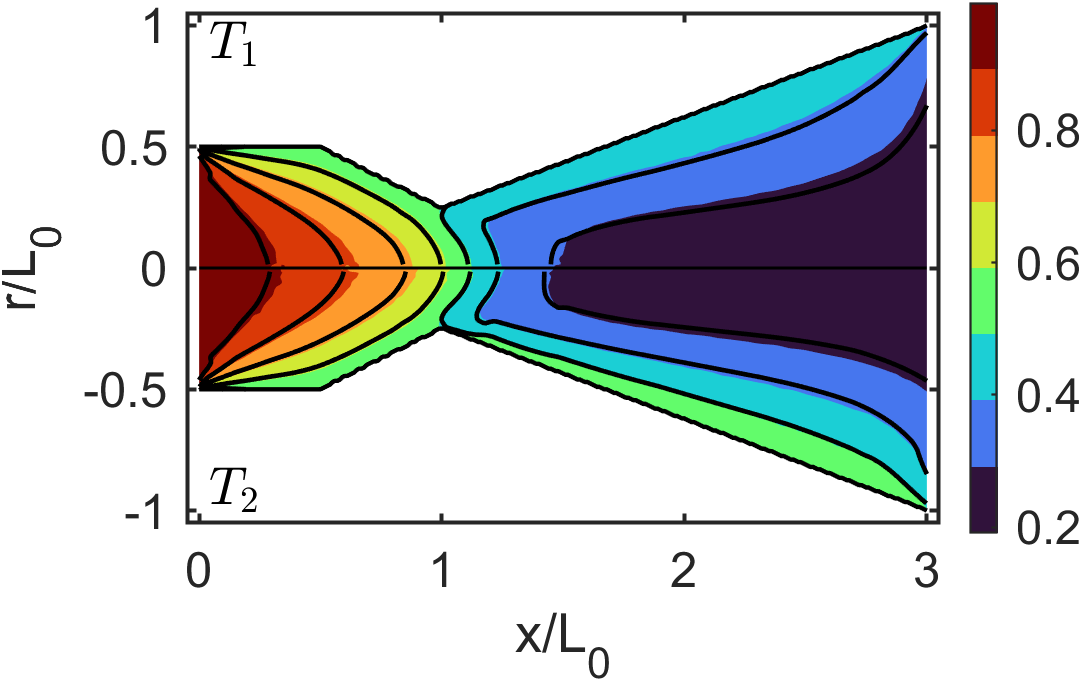}
\end{minipage}
\caption{
Species number density $n_s/n_0$ (top), axial velocity $u_{s,x}/v_0$ (middle), and temperature $T_s/T_0$ (bottom) in the nozzle meridional plane at $\Kn_1=0.1$ (left) and 0.01 (right). Filled contours denote DSMC results and black isolines denote GSIS-DG. The upper and lower halves of each panel correspond to species 1 and 2, respectively.
}
\label{fig:dg-nozzle-fields}
\end{figure}

Figure~\ref{fig:dg-nozzle-fields} compares the species number density, axial velocity, and temperature obtained by GSIS-DG and DSMC in the meridional plane. Background contours denote DSMC results, while black isolines denote GSIS-DG; the upper and lower halves of each panel correspond to the light and heavy species, respectively. Both species accelerate through the contraction and downstream expansion, accompanied by a substantial decrease in number density and temperature. The light species reaches a higher downstream axial velocity, particularly at $\Kn_1=0.1$, whereas the velocity difference between the two species is reduced at $\Kn_1=0.01$ because of stronger interspecies momentum exchange. Overall, GSIS-DG reproduces the principal DSMC flow structures on the very coarse three-dimensional physical mesh. The most visible deviations occur in the heavy-species density near the throat.

\begin{figure}[!htbp]
\centering
\makebox[\textwidth][c]{%
\includegraphics[width=0.4\textwidth,keepaspectratio]{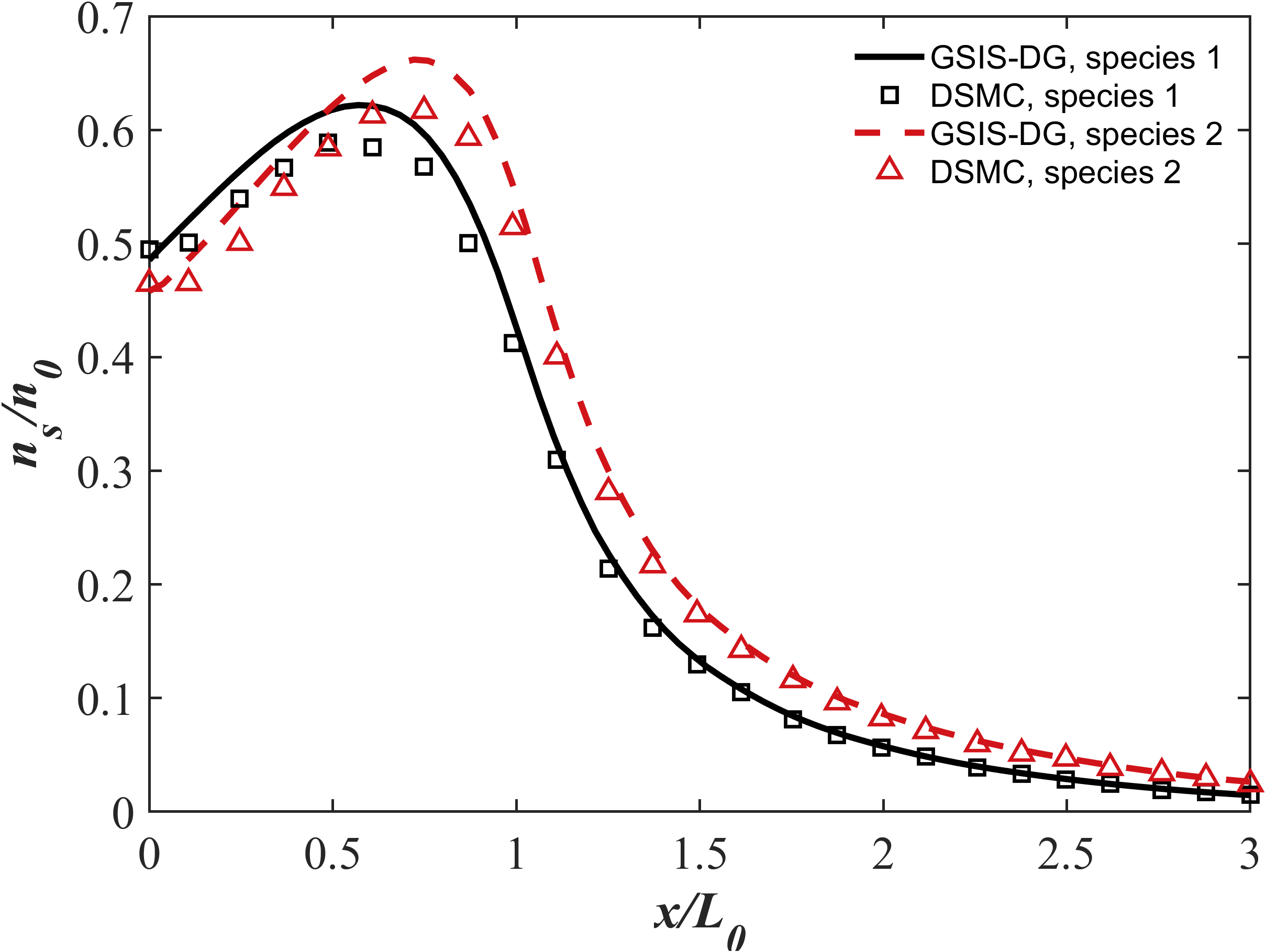}%
\includegraphics[width=0.4\textwidth,keepaspectratio]{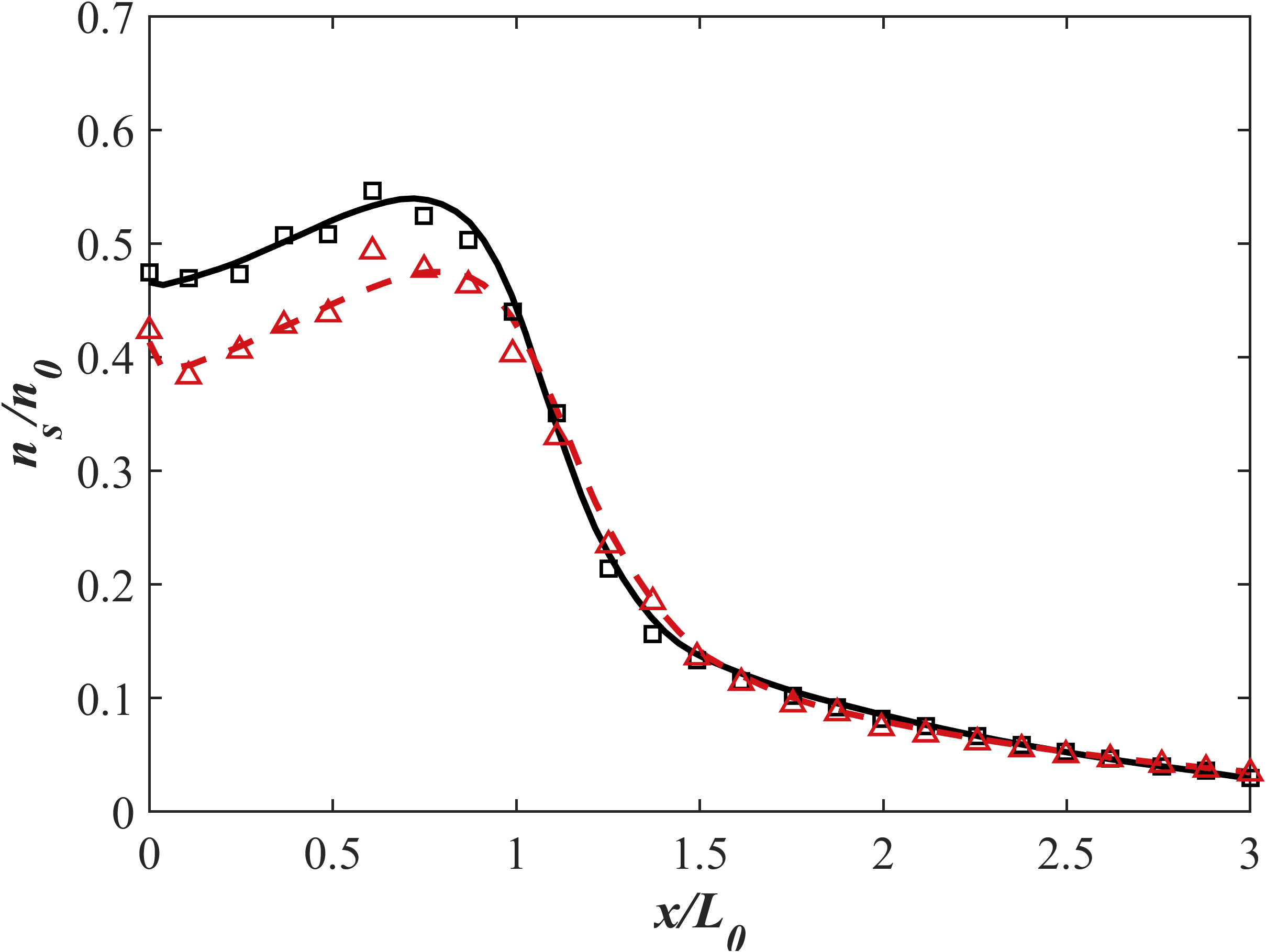}%
}\\[0.25em]
\makebox[\textwidth][c]{%
\includegraphics[width=0.4\textwidth,keepaspectratio]{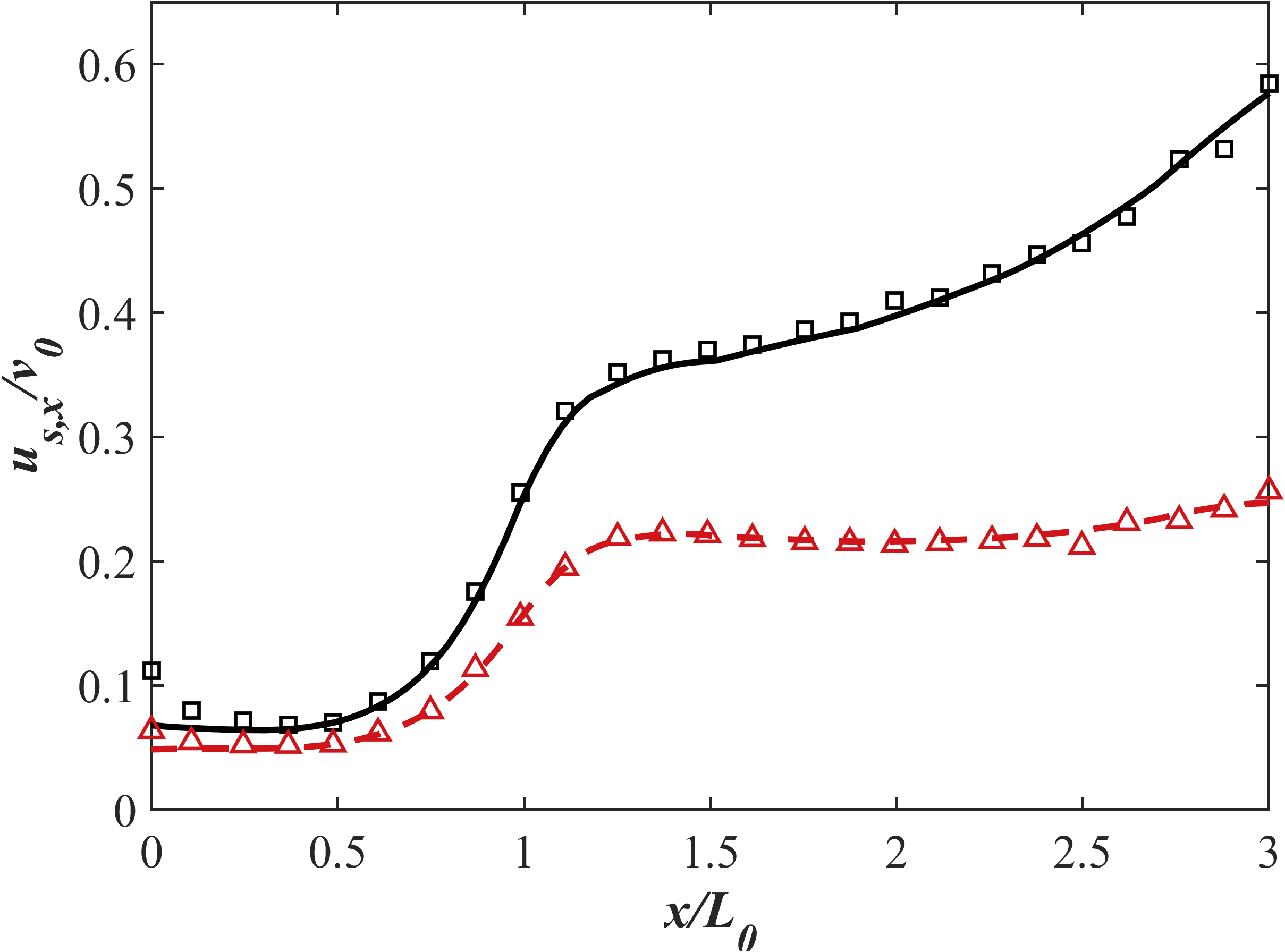}%
\includegraphics[width=0.4\textwidth,keepaspectratio]{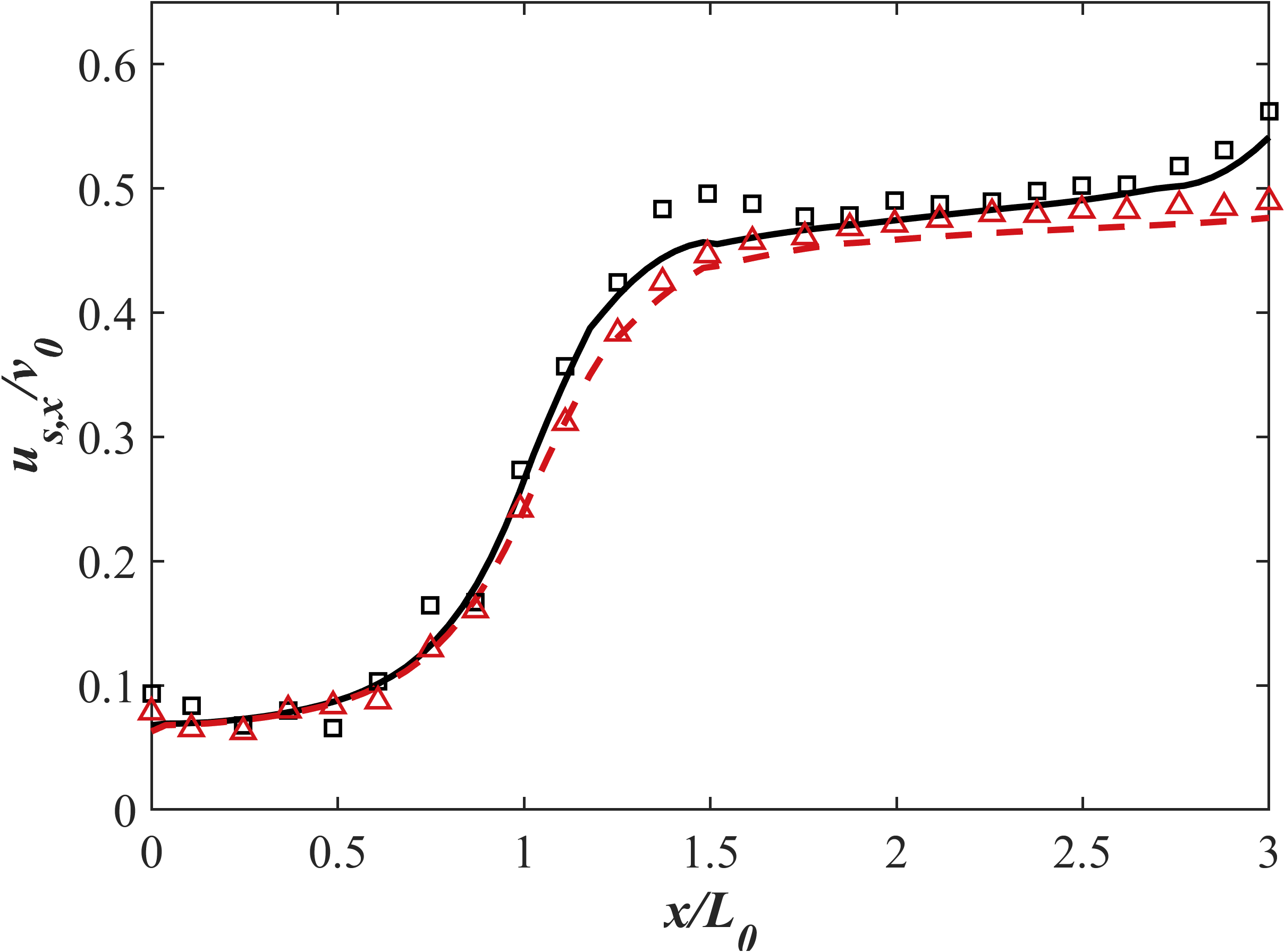}%
}
\caption{
Species number densities (upper row) and axial velocities (lower row) along the nozzle centerline at $\Kn_1=0.1$ (left) and 0.01 (right). Lines denote GSIS-DG and symbols denote DSMC; black and red correspond to the light and heavy species, respectively.
}
\label{fig:dg-nozzle-centerlines}
\end{figure}

The centerline profiles in Fig.~\ref{fig:dg-nozzle-centerlines} provide a more quantitative comparison. The species number densities first increase in the converging section and then decrease by more than one order of magnitude during the downstream expansion. At $\Kn_1=0.1$, the density maximum of the heavy species is shifted slightly downstream relative to that of the light species. The light species continues to accelerate downstream of the throat, whereas the heavy-species velocity varies more gradually. At $\Kn_1=0.01$, the two axial velocities become considerably closer, consistent with the stronger interspecies coupling. The GSIS-DG predictions agree well with the DSMC data over most of the nozzle. The more noticeable deviations in the heavy-species density appear near the throat, around $x/L_0\simeq0.8$ for $\Kn_1=0.1$ and $x/L_0\simeq1.2$ for $\Kn_1=0.01$. These discrepancies may be partly associated with differences between the collision descriptions in the kinetic model and DSMC.

\begin{figure}[!htbp]
\centering
\makebox[\textwidth][c]{%
\includegraphics[width=0.50\textwidth,keepaspectratio]{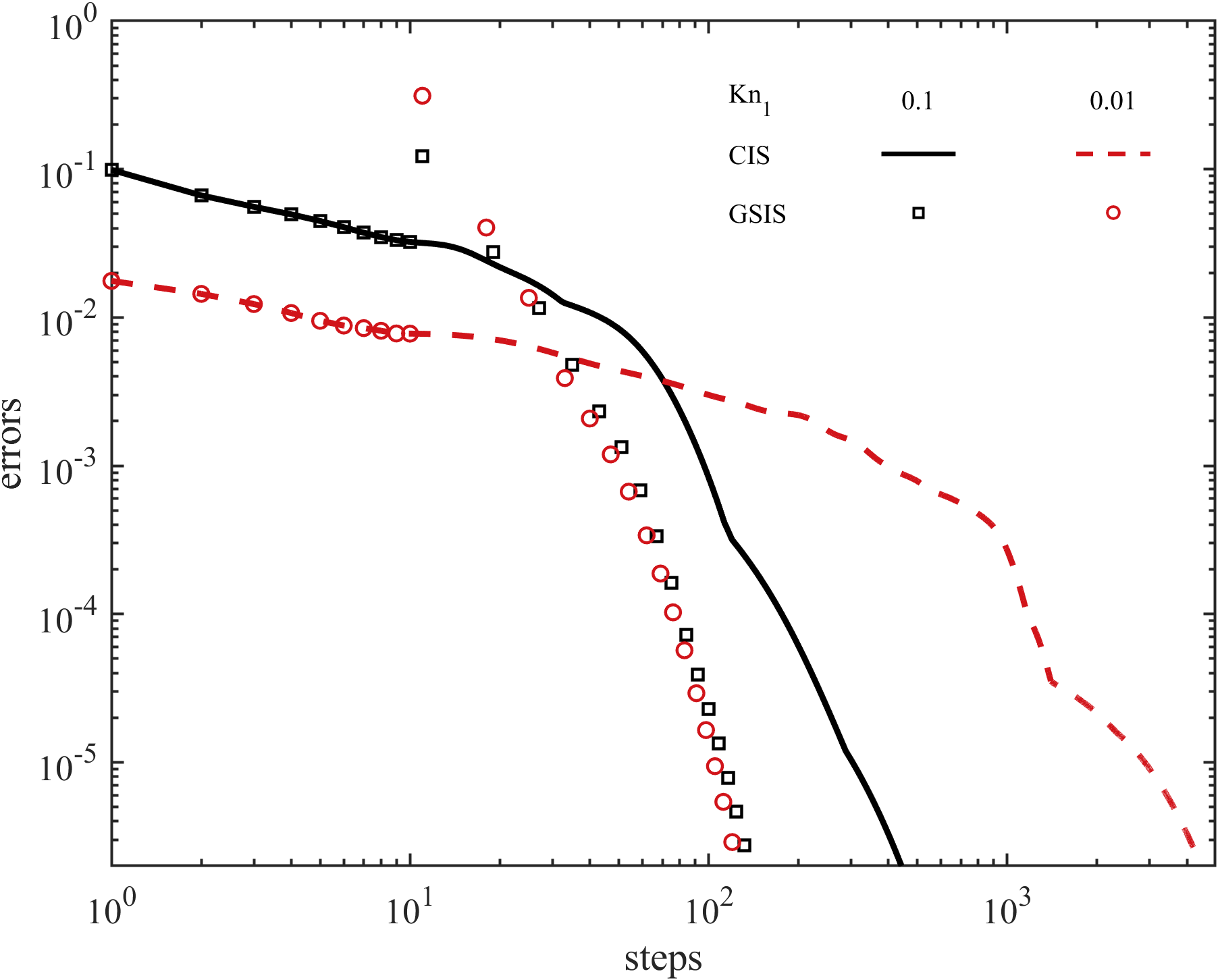}
}
\caption{
Convergence histories of CIS-DG and GSIS-DG for the nozzle flow at $\Kn_1=0.1$ and 0.01. Black and red denote $\Kn_1=0.1$ and 0.01, respectively; solid and dashed lines show CIS-DG for these two cases, and matching open squares for $\Kn_1=0.1$ and open circles for $\Kn_1=0.01$ show GSIS-DG. The horizontal axis counts total DVM updates, including the ten initial CIS updates in each GSIS-DG calculation.
}
\label{fig:dg-nozzle-convergence}
\end{figure}

The convergence histories in Fig.~\ref{fig:dg-nozzle-convergence} and the computational costs in Table~\ref{tab:dg-nozzle-cost} demonstrate that the rapid convergence of GSIS-DG is retained in the three-dimensional calculation. The horizontal axis in Fig.~\ref{fig:dg-nozzle-convergence} counts kinetic updates, including the ten initial CIS updates used to initialize GSIS-DG. At $\Kn_1=0.1$, CIS-DG and GSIS-DG require 502 and 148 kinetic updates, respectively, while the corresponding computational costs are 12.8 and 6.2 core hours. Thus, GSIS-DG reduces the iteration count by a factor of 3.4 and the computational cost by a factor of 2.1. At $\Kn_1=0.01$, the difference becomes much more pronounced. GSIS-DG converges within 134 kinetic updates and 5.9 core hours, whereas CIS-DG remains unconverged after 5000 updates and 129.2 core hours. Therefore, the speedup exceeds 37 in kinetic updates and 21 in computational cost. The smaller speedup in core hours than in update count reflects the additional solution of the macroscopic synthetic equations in each GSIS-DG outer iteration. Nevertheless, the results demonstrate that the synthetic iteration effectively removes the severe convergence deterioration of conventional kinetic iterations as the flow approaches the near-continuum regime, even for a fully three-dimensional nozzle flow.

For reference, the axisymmetric DSMC simulations reported in Ref.~\cite{zeng2024gsismixture} consumed approximately 5 and 1030 core hours at $\Kn_1=0.1$ and 0.01, respectively. At $\Kn_1=0.1$, the fully three-dimensional GSIS-DG calculation requires 6.2 core hours, which is comparable to the reported cost of the two-dimensional axisymmetric DSMC calculation. At $\Kn_1=0.01$, the reported GSIS-DG cost is more than two orders of magnitude lower than the literature-reported DSMC cost, despite the present calculation being fully three-dimensional. This comparison is intended only as an indication of the computational scale, since the DSMC data are taken from the literature and were obtained with different hardware, dimensionality, and implementation details.


\begin{table}[!htbp]
\begin{threeparttable}
\centering
\caption{
Computational performance of CIS-DG and GSIS-DG for the nozzle flow. The physical domain is discretized by 429 hexahedral elements with $p=2$. All present calculations use 40 MPI ranks on a dual-socket node equipped with two 32-core Intel Xeon 7285H processors at 2.5~GHz. The total computational cost is reported in core hours, and the GSIS-DG update counts include ten initial CIS updates.
}
\label{tab:dg-nozzle-cost}
\scriptsize
\begin{tabular*}{0.94\textwidth}{@{\extracolsep{\fill}}c c cc cc cc@{}}
\toprule
& DSMC\tnote{1} & \multicolumn{2}{c}{CIS-DG\tnote{2}} & \multicolumn{2}{c}{GSIS-DG} & \multicolumn{2}{c}{Speedup ratio} \\
\cmidrule(lr){3-4}\cmidrule(lr){5-6}\cmidrule(lr){7-8}
$\Kn_1$ & cost & steps & cost & steps & cost & in steps & in cost \\
\midrule
0.1  & 5    & 502       & 12.8     & 148 & 6.2 & 3.4     & 2.1 \\
0.01 & 1030 & $>5000$    & $>129.2$   & 134 & 5.9 & $>37.3$   & $>21.9$ \\
\bottomrule
\end{tabular*}
\begin{tablenotes}
\footnotesize
\item[1] The DSMC costs are taken from Ref.~\cite{zeng2024gsismixture}, where two-dimensional axisymmetric simulations were conducted. These values are literature references, not timings measured in the present comparison. At $\Kn_1=0.1$, 6165 cells and 73,000 particles were employed; at $\Kn_1=0.01$, 373,761 cells and 2,962,500 particles were employed.
\item[2] For $\Kn_1=0.01$, CIS remains unconverged at 5000 updates. The entries with $>$ give conservative lower bounds for convergence to the prescribed tolerance.
\end{tablenotes}
\end{threeparttable}
\end{table}

\section{Conclusion}\label{sec:5}

The present results demonstrate that high-order DG discretization improves the spatial efficiency of GSIS for rarefied gas mixture flows. The asymptotic-preserving analysis establishes conditions under which the multicomponent Navier-Stokes limit is recovered on meshes that do not resolve the molecular mean free path, while the numerical results show that increasing the polynomial degree reduces the number of elements needed to resolve the flow structures accurately. Comparisons with finite-volume GSIS further show that the reduction in element count can translate into fewer total spatial degrees of freedom and lower computational cost. Meanwhile, the synthetic equations maintain rapid convergence toward the steady state, with an increasing advantage over conventional kinetic iteration as the Knudsen number decreases. Together, these findings demonstrate that the spatial accuracy can be improved without sacrificing the asymptotic-preserving and fast-converging properties of GSIS.

The present method provides an efficient approach for resolving species dependent non-equilibrium behavior in rarefied gas mixture flows. The nozzle calculations demonstrate that differences in species density, velocity, and temperature can be captured on coarse physical meshes, while rapid convergence is maintained toward the continuum regime. These capabilities suggest potential applications to larger scale mixture flows in which both spatial resolution and convergence efficiency are demanding. We believe that the present framework can be further extended, through $hp$ adaptation and suitable stabilization, to tackle more challenging high-Mach-number mixture flows.

\section*{Acknowledgments}

The work is financially supported by the Hong Kong Research Grants Council under the Early Career Scheme (project No. 26211523).

\appendix

\bibliographystyle{elsarticle-num}
\bibliography{ref}
\end{document}